\documentclass[twoside,11pt,reqno]{amsart}
\usepackage{amsmath,amssymb,amscd,mathrsfs,epic,wasysym,latexsym,tikz,mathrsfs,cite,hyperref,enumerate}
\usepackage{pb-diagram}
\usepackage[all]{xy}
\usepackage{mathdots}

\makeatletter

\numberwithin{equation}{section}
\allowdisplaybreaks

\newtheorem{Proposition}[equation]{Proposition}
\newtheorem{Lemma}[equation]{Lemma}
\newtheorem{Theorem}[equation]{Theorem}
\newtheorem{Corollary}[equation]{Corollary}
\newtheorem{MainTheorem}{Theorem}

\theoremstyle{definition}  %% makes all of the theorem environments which follow appear in \rm

\newtheorem{Remark}[equation]{Remark}

\theoremstyle{remark}

\let\<\langle
\let\>\rangle

\newcommand\Comment[2][\relax]{\space\par\medskip\noindent%
   \fbox{\begin{minipage}{\textwidth}\textbf{Comment\ifx\relax#1\else---#1\fi}\newline%
        #2\end{minipage}}\medskip
}

\def\bi{\text{\boldmath$i$}}

\def\ba{\text{\boldmath$a$}}
\def\bb{\text{\boldmath$b$}}

\def\pmod#1{\text{ }(\text{\rm mod } #1)\,}

\newcommand{\Hom}{\operatorname{Hom}}
\newcommand{\Ext}{\operatorname{Ext}}

\newcommand{\End}{\operatorname{End}}

\newcommand{\cont}{\operatorname{cont}}

\def\sgn{\mathtt{sign}}

\newcommand{\res}{\operatorname{res}}

\newcommand{\infl}{\operatorname{infl}}

\newcommand{\Z}{\mathbb{Z}}

\newcommand{\0}{{\bar 0}}
\renewcommand{\1}{{\bar 1}}
\renewcommand{\2}{{\bar 2}}
\def\eps{{\varepsilon}}
\def\phi{{\varphi}}

\newcommand{\F}{{\mathbb F}}

\newcommand{\rank}{{\mathfrak r}}

\newcommand{\level}{{\mathfrak l}}

\newcommand{\type}{{\mathfrak t}}

\newcommand{\ga}{\gamma}

\newcommand{\la}{\lambda}
\newcommand{\La}{\Lambda}
\newcommand{\al}{\alpha}
\newcommand{\be}{\beta}

\def\Si{\mathfrak{S}}

\newcommand{\de}{\delta}
\newcommand{\De}{\Delta}

\def\triv#1{\O_{#1}}

\newcommand{\GL}{{\sf GL}}
\newcommand{\T}{{\sf T}}

\newcommand{\Irr}{{\mathrm {Irr}}}

\newcommand{\Ind}{{\mathrm {Ind}}}

\newcommand{\diag}{{\mathrm {diag}}}
\newcommand{\Res}{{\mathrm {Res}}}

\newcommand{\C}{{\mathbb C}}

\newcommand{\dar}{{\downarrow}}

\renewcommand{\mod}{\bmod \,}

\newcommand{\AAA}{\mathfrak{A}}

\def\Parp{{\mathscr P}_p}
\def\Par{{\mathscr P}}

\def\k{\Bbbk}

\def\T{{\mathtt T}}

\def\op{{\mathrm{op}}}

\def\onto{{\twoheadrightarrow}}
\def\into{{\hookrightarrow}}

\def\mod#1{#1\!\operatorname{-mod}}

\def\iso{\stackrel{\sim}{\longrightarrow}}

\def\Pol{{\mathcal P}}

\def\triv{{\tt triv}}
\def\col{{\tt col}}
\def\row{{\tt row}}
\def\ch{\operatorname{ch}}

\def\H{\mathsf{H}}
\def\s{\mathsf{s}}
\def\x{\mathsf{x}}
\newcommand{\trinom}[4]{\left(\begin{array}{@{}c@{}} #1\\ #2, #3, #4\end{array}\right)}

{\catcode`\|=\active
  \gdef\set#1{\mathinner{\lbrace\,{\mathcode`\|"8000%
  \let|\midvert #1}\,\rbrace}}
}
\def\midvert{\egroup\mid\bgroup}

\colorlet{darkgreen}{green!50!black}
\tikzset{dots/.style={very thick,loosely dotted},
         greendot/.style={fill,circle,color=darkgreen,inner sep=1.5pt,outer sep=0},
         blackdot/.style={fill,circle,color=black,inner sep=1.5pt,outer sep=0},
         graydot/.style={fill,circle,color=gray,inner sep=1.1pt,outer sep=0}
}
\def\greendot(#1,#2){\node[greendot] at(#1,#2){}}
\def\blackdot(#1,#2){\node[blackdot] at(#1,#2){}}
\def\graydot(#1,#2){\node[graydot] at(#1,#2){}}

\newenvironment{braid}{% sets defaults for the braid diagrams
  \begin{tikzpicture}[baseline=6mm,black,line width=1pt, scale=0.32,
                      draw/.append style={rounded corners},
                      every node/.append style={font=\fontsize{5}{5}\selectfont}]%
  }{\end{tikzpicture}
}

\def\Grid(#1,#2){%  draws a coordinate grid inside a braid diagram
  \draw[very thin,gray,step=2mm] (0,0)grid(#1,#2);
  \draw[very thin,darkgreen,step=10mm] (0,0)grid(#1,#2);
}

\newcommand\Tableau[2][\relax]{
  \begin{tikzpicture}[scale=0.5,draw/.append style={thick,black}]
    \ifx\relax#1\relax%
    \else % shade the boxes in #1
      \foreach\box in {#1} { \filldraw[blue!30]\box+(-.5,-.5)rectangle++(.5,.5); }
    \fi
    \newcount\row\newcount\col
    \row=0
    \foreach \Row in {#2} {
       \col=1
       \foreach\k in \Row {
          \draw(\the\col,\the\row)+(-.5,-.5)rectangle++(.5,.5);
          \draw(\the\col,\the\row)node{\k};
          \global\advance\col by 1
       }
       \global\advance\row by -1
    }
  \end{tikzpicture}
}

\newcommand\YoungDiagram[2][\relax]{
  \begin{tikzpicture}[scale=0.5,draw/.append style={thick,black}]
    \ifx\relax#1\relax%
    \else % shade the boxes in #1
    \foreach\box in {#1} {
      \filldraw[blue!30]\box rectangle ++(1,1);
    }
    \fi
    \newcount\row
    \row=0
    \foreach \col in {#2} {
       \draw(1,\the\row)grid ++(\col,1);
       \global\advance\row by -1
    }
  \end{tikzpicture}
}

\newdimen\hoogte    \hoogte=12pt    
\newdimen\breedte   \breedte=14pt  
\newdimen\dikte     \dikte=0.5pt 

\newenvironment{Young}{\begingroup
       \def\vr{\vrule height0.89\hoogte width\dikte depth 0.2\hoogte}
       \def\fbox##1{\vbox{\offinterlineskip
                    \hrule height\dikte
                    \hbox to \breedte{\vr\hfill##1\hfill\vr}
                    \hrule height\dikte}}
       \vbox\bgroup \offinterlineskip \tabskip=-\dikte \lineskip=-\dikte
            \halign\bgroup &\fbox{##\unskip}\unskip  \crcr }
       {\egroup\egroup\endgroup}
\def\Youngdiagram#1{\relax\ifmmode\vcenter{\,\begin{Young}#1\end{Young}\,}\else%
              $\vcenter{\,\begin{Young}#1\end{Young}\,}$\fi}

\begin{document}

\title[Level, rank and tensor growth]{{\bf Level, rank, and tensor growth of representations of symmetric groups}}

\author{\sc Alexander Kleshchev}
\address{Department of Mathematics\\ University of Oregon\\Eugene\\ OR~97403\\ USA}
\email{klesh@uoregon.edu}

\author{\sc Michael Larsen}
\address
{Department of Mathematics\\ Indiana University\\ Bloomington\\ IN~47405\\ USA} 
\email{mjlarsen@indiana.edu}

\author{\sc Pham Huu Tiep}
\address
{Department of Mathematics\\ Rutgers University\\ Piscataway\\ NJ~08854\\ USA} 
\email{tiep@math.rutgers.edu}

%\subjclass[2010]{20C30, 20C20}

\thanks{The first author was supported by the NSF grant DMS-2101791, Charles Simonyi Endowment at the Institute for Advanced Study, and the Simons Foundation.
The second author was partially supported by the NSF grant DMS-2001349 and the Simons Foundation.
The third author was partially supported by the NSF (grants DMS-1840702 and DMS-2200850), the Simons Foundation, and 
the Joshua Barlaz Chair in Mathematics.}

\thanks{Part of this work was done when the authors visited the Institute for Advanced Study and Princeton University. 
It is a pleasure to thank the Institute for Advanced Study and Princeton University for generous hospitality and stimulating 
environment.} 

\begin{abstract}
We develop a theory of levels for irreducible representations of symmetric groups $\Si_n$ analogous to the theory of levels for finite classical groups. 
A key property of level is that the level of a character, provided it is not too big compared to $n$, gives a good lower bound on its degree, and, moreover, every character of low degree is either itself of low level or becomes so after tensoring with the sign character. Furthermore, if $l_1$ and $l_2$ satisfy a linear upper bound in $n$, then the maximal level of composition factors of the 
tensor product of representations of levels $l_1$ and $l_2$ is 
$l_1+l_2$. To prove all of this in positive characteristic, we develop the notion of rank, which is an analogue of the notion of rank of cross-characteristic 
representations of finite classical groups. We show, using modular branching rules and degenerate affine Hecke algebras, 
that the level and the rank 
agree, as long as the level is not too large. We exploit Schur-Weyl  duality, modular Littlewood-Richardson coefficients and tilting modules to prove a modular analogue of the  Murnaghan-Littlewood theorem on Kronecker products for symmetric groups. %A by-product of the rank approach gives the new dimension bound for modular representations of $\Si_n$. 
As an application, we obtain representation growth results for both ordinary and modular representations of symmetric and alternating groups analogous to those for finite groups of Lie type.  
\end{abstract}

\maketitle

\tableofcontents

\section{Introduction}
In this paper, we develop a theory of {\it levels} for irreducible representations of symmetric groups $\Si_n$,
analogous to the theory for finite classical groups developed in \cite{GLT2, GLT3}.

For a partition $\la=(\la_1,\la_2,\dots)$ of $n$ we define the {\em level} of $\la$ to be 
$$
\level(\la):=n-\la_1.
$$

The irreducible $\C\Si_n$-modules $S^\lambda$ are indexed by partitions $\lambda$ of $n$. We refer to $\level(\la)$ as the level of $S^\la$ or the level of the corresponding irreducible character $\chi^\la$. 
For fixed level $l$ and varying $n$, the Hook Length Formula shows that  the dimension of an irreducible $\C\Si_n$-module of level $l$ is a polynomial in $n$ of degree~$l$.  %These $S^\lambda$ have {\it level\, $l$}.  

A key property of level, from our point of view, is that the level of a character, provided it is not too big compared to $n$, gives a good lower bound on its degree, and, moreover, every character of low degree is either itself of low level or becomes so after tensoring with the sign character. Furthermore, if $l_1$ and $l_2$ satisfy a linear upper bound in $n$, then the maximal level of composition factors of the 
tensor product of representations of levels $l_1$ and $l_2$\, is\, 
$l_1+l_2$. We refer to this last property as the {\em additivity property of levels}
%Furthermore, if $k_1$ and $k_2$ satisfy a linear upper bound in $n$, the tensor product of representations of levels $k_1$ and $k_2$ decomposes as a direct sum of representations of level $\le k_1+k_2$.  

In positive characteristic, we can still define the level of an irreducible representation of $\Si_n$.  We can still obtain lower bounds (somewhat weaker than in the characteristic $0$ case) for the degree of a representation of specified level.  Moreover, the additivity property of levels still holds. 
%if $k_1$ and $k_2$ satisfy a linear upper bound in $n$, then the maximal level of composition factors of the tensor product of representations of levels $k_1$ and $k_2$ is $k_1+k_2$. 
To prove all of this, we develop the notion of {\it rank}, which is an analogue of the notion of rank of cross-characteristic 
representations of finite classical groups \cite{GLT2,GLT3}. 

Let $\F$ be a field of arbitrary characteristic $p\geq 0$, and in this introduction let us assume for simplicity that $\F$ is algebraically closed (this is not needed for most results of the paper). Let $V$ be an $\F\Si_n$-module. 
If $p\neq 2$, we define the {\em $2$-rank} of  $V$ as the largest  $r$ such that there exist distinct $i_1,j_1,\dots,i_r,j_r\in\{1,\dots,n\}$ and a simultaneous $(-1)$-eigenvector in $V$ for the transpositions $(i_1,j_1),\dots,(i_r,j_r)$.  
The $2$-rank of $V$ will be denoted $\rank_2(V)$. If $p\neq 3$, let $\zeta$ be a primitive $3$rd root of unity in $\F$, and define the {\em $3$-rank} of $V$ as the largest $r$ such that there exist distinct $i_1,j_1,k_1,\dots,i_r,j_r,k_r\in\{1,\dots,n\}$ and a simultaneous  $\zeta$-eigenvector in $V$ for the $3$-cycles $(i_1,j_1,k_1),\dots,(i_r,j_r,k_r)$. The $3$-rank of $V$ will be denoted $\rank_3(V)$.

The irreducible $\F\Si_n$-modules $D^\lambda$ are indexed by $p$-regular partitions $\lambda$ of $n$ (if $p=0$, all partitions are $p$-regular by definition and $D^\la=S^\la$). %We again refer to $\level(\la)$ as the level of $D^\la$. 
Our first main result (see Theorems~\ref{TRL} and \ref{TRL2})  shows in particular that the level and the rank 
agree as long as the level is not too large:

\begin{MainTheorem} \label{TA}
Let $\la$ be a $p$-regular partition of $n$. 
\begin{enumerate}
\item[{\rm (i)}] If $p\neq 2$ then 
$\rank_2(D^\la)=\min\{\level(\la),\lfloor n/2\rfloor\}.$ 
\item[{\rm (ii)}] If $p= 2$ then 
$\rank_3(D^\la)=\min\{\level(\la),\lfloor n/3\rfloor\}.$ 
\end{enumerate}
\end{MainTheorem}

The proof of Theorem~\ref{TA} uses modular branching rules \cite{KBrII,KDec} and some basic representation theory of degenerate affine Hecke algebras. 

Our second key result (see 
Theorem \ref{TMur}) is a generalization to an arbitrary characteristic of the classical Murnaghan-Littlewood theorem, which describes some important (modular) Kronecker coefficients in terms of (modular) Littlewood-Richardson coefficients. For a partition $\la=(\la_1,\la_2,\dots,\la_h)$ we define the partition 
$$
 \bar\la:=(\la_2,\dots,\la_h)
$$
of $\level(\la)$. 

\begin{MainTheorem} \label{TB}
Let $\la,\mu,\nu$ be $p$-regular partitions of $n$,\, $l=\level(\la)$\, and\, $m=\level(\mu)$. If\, $\level(\nu)=l+m$ then 
 $$
 [D^\la\otimes D^\mu: D^\nu]=[\Ind_{\Si_{l}\times \Si_{m}}^{\Si_{l+m}} (D^{\bar\la}\boxtimes D^{\bar\mu}): D^{\bar\nu}].
 $$
\end{MainTheorem}

Theorem~\ref{TB} substantially strengthens the additivity property of levels. Its proof uses Schur-Weyl duality, equality of some modular Littlewood-Richardson coefficients established in \cite{BKLR}, and the theory of tilting modules. We believe this theorem will have further applications. 

As a by-product of the rank approach, 
we obtain a new dimension bound for modular representations of $\Si_n$ (see Theorem \ref{TBound}) :

\begin{MainTheorem}\label{TC}
Let $\la$ be a $p$-regular partition of $n$ of level $l$. Then
$$
\dim D^\la \geq
\left\{
\begin{array}{ll}
\displaystyle\binom{\lfloor n/2 \rfloor}{l} &\hbox{if $p\neq 2$,}
\vspace{2mm}
\\
\displaystyle2^r\binom{\lfloor n/3 \rfloor}{l} &\hbox{if $p=2$.}
\end{array}
\right.
$$
\end{MainTheorem}

 In some crucial situations Theorem~\ref{TC} yields stronger results than lower bounds coming from \cite{KMT}, see \S\ref{SSBound} for further discussion of this. In particular, the new bound is needed for the applications to representation growth described below.

As a first application of the results on rank and level described above, we establish some representation growth results for both ordinary and modular representations of symmetric and alternating groups analogous to those proved for finite groups of Lie type in \cite{LST}.  
This means that if $V$ and $W$ are representations of a group $G$ of Plancherel measures $|V|$ and $|W|$ small compared to $|G|$, then the Plancherel measure $|V\otimes W|$ is large compared to 
$\sqrt{|V|\, |W|}$, see \S5 for more details. 
%However, these results are not as strong as the analogous results for groups of Lie type.
In characteristic zero, our bounds apply when the level of every irreducible constituent of $V$ and $W$ satisfies a linear upper bound in $n$,
so, in particular, when $|V|$ and $|W|$ are bounded above by an exponential in $n$ with base close enough to $1$. In positive characteristic, our results hold only if the levels of the constituents satisfy a power bound $n^\alpha$ for 
some $0 < \alpha < 1$.

\begin{MainTheorem} \label{TD}
Let $\epsilon > 0$ and $G_n$ be the symmetric group $\Si_n$ or the alternating group $\AAA_n$. 
\begin{enumerate}
\item[{\rm (i)}] there exists $\delta=\delta(\epsilon) > 0$ and $N=N(\epsilon)$ such that for all $n\geq N$ and $\C G_n$-modules $V,W$  with $|V|,|W| < (1+\delta)^n$, we have  
$$|V\otimes W| \ge \big(|V|\,|W|\big)^{1-\epsilon}.$$ 
In particular, $|V \otimes V| \geq |V|^{2-2\eps}$.

\item[{\rm (ii)}] there exists $N=N(\eps)$ such that for all $n\geq N$ and $\F G_n$-modules\, $V,W$  with $|V|, |W| < 2^{n^{2/3-\eps}}$, we have 
$$|V\otimes W| \ge \bigl(|V| \, |W|\bigr)^{\frac{1}{2}+\frac{\eps}{6}}.$$ 
In particular, $|V \otimes V| \geq |V|^{1+\eps/3}$.

%\item[{\rm (ii)}] For all $k \in \Z_{\geq 2}$ there exists $N=N(\eps,k) > 1$ such that for all $n \geq N$ and $\F G_n$-modules $V_1, \ldots, V_k$ with $|V_1|,\dots, |V_k| < 2^{n^{(2k-2)/(2k-1)-\eps}}$ we have  $$|V_1 \otimes V_2 \otimes \ldots \otimes V_k| \ge \bigl(|V_1| \, |V_2| \, \cdots \, |V_k|\bigr)^{(1+\eps/3)/k}.$$ In particular, $|V^{\otimes k}| \geq |V|^{1+\eps/3}$.

\item[{\rm (iii)}] Assuming\, $\eps <1$, there exists $N=N(\eps)$ such that for all $n \geq N$ and $\F G_n$-modules $V,W$ 
 with $|V|, |W| < 2^{n^{\eps/2}}$, we have 
$$|V \otimes W| \ge \bigl(|V| \, |W|\bigr)^{1-\eps}.$$
In particular, $|V \otimes V| \geq |V|^{2-2\eps}$.
\end{enumerate} 
\end{MainTheorem}

Theorem~\ref{TD} is proved in Theorems~\ref{main1},\,\ref{main-2sa}, \ref{alt} and Corollary~\ref{main-2s}. See also Theorem~\ref{main-3s} and Remark~\ref{RKStep} for a generalization to more than two tensor factors.

\section{Preliminaries}
\subsection{General notation}
Throughout the paper, $\F$ denotes an arbitrary field of arbitrary characteristic $p\geq 0$. If $p>0$, we set
$$
I:=\{0,1,\dots,p-1\}
$$
identified with $\Z/p\Z$ which is in turn identified with the prime subfield of $\F$. In particular every $i\in I$ makes sense as an element of $\F$. We use the total order  $0<1<\dots<p-1$ on~$I$. 
If $p=0$, we interpret $I$ as $\Z$ embedded naturally into $\F$. 

All modules in this paper are assumed to be finite-dimensional.  
Let $A$ be an $\F$-algebra. We denote by $\mod{A}$ the category of finite-dimensional $A$-modules. If $V,L,T\in\mod{A}$ with $L$ irreducible and $T$ indecomposable, we denote by $[V:L]$ the multiplicity of $L$ as a composition factor of $V$ (well-defined by Jordan-H\"older), and by $(V:T)$ the multiplicity of $T$ as a direct summand of $V$ (well-defined by Krull-Schmidt). 

Let $B$ be another $\F$-algebra. Given $V\in\mod{A}$ and $W\in\mod{B}$, we have the outer tensor product module $V\boxtimes W\in\mod{A\otimes B}$. If $B$ is a subalgebra of $A$, we have the functors of restriction and induction
$$
\Res^A_B:\mod{A}\to\mod{B}\quad \text{and}\quad \Ind_B^A:\mod{B}\to\mod{A}.
$$

Let $G$ be a group. We denote by $\triv_G$, or simply $\triv$ when $G$ is clear, the trivial $\F G$-module. 
If $H$ is a subgroup of $G$ and $V\in\mod{\F G}$, we often write $V\dar_H$ instead of $\Res^{\F G}_{\F H}$ and $\Ind_H^G$ for $\Ind^{\F G}_{\F H}$. 

We denote by $\Si_n$ the symmetric group on $n$ letters with transpositions $(i,j)$, $3$-cyles $(i,j,k)$, etc. We interpret $\Si_0$ as the trivial group. 
%We also denote  $s_r:=(r,r+1)$ for $1\leq r<n$. 
The {\em $i$th Murphy element} is 
$$L_i:=\sum_{j=1}^{i-1}(j,i)\in\F\Si_n.$$
For $k_1+\dots+k_l=n$, the group $\Si_{k_1}\times\dots\times \Si_{k_l}$ is always considered as a standard Young subgroups of $\Si_n$.
We denote by $\sgn_{\Si_n}$, or simply $\sgn$ when $n$ is clear, the sign representation of $\F \Si_n$. For a subgroup $G\leq \Si_n$, we still have the representation $\sgn:=\sgn_{\Si_n}\dar_G$.  

We denote by $\Par(n)$ the set of all {\em partitions} of $n$ and by $\Parp(n)$ the set of all {\em $p$-regular} partitions of $n$, see \cite[10.1]{JamesBook}. The conjugate partition of $\la$ is denoted $\la'$, see \cite[Definition~3.5]{JamesBook}. 
We have a {\em dominance order} $\unrhd$ on partitions, see \cite[3.2]{JamesBook}. Given $\la\in\Par(l)$ and $\mu\in\Par(m)$, we have the {\em sum partition} $\la+\mu\in\Par(l+m)$, which is defined as the partition $\nu$ with $\nu_i=\la_i+\mu_i$ for all $i$. 

We have the family 
$
\{S^\la\mid\la\in\Par(n)\}
$
of  Specht modules over $\F\Si_n$, see \cite[\S4]{JamesBook}, and the family $
\{D^\la\mid\la\in\Parp(n)\}
$
of irreducible modules over $\F\Si_n$, see \cite[\S11]{JamesBook}. If $p=0$ we have $D^\la=S^\la$ for all $\la\in \Par_0(n)=\Par(n)$.

\subsection{Good and normal nodes}
We identify a partition $\la=(\la_1,\la_2,\dots)$ with its {\em Young diagram} 
$
\{(r,s)\in\Z_{>0}\times\Z_{>0}\mid s\leq\la_r\}.
$
The elements $(r,s)\in\Z_{>0}\times\Z_{>0}$ are referred to as {\em nodes}. We refer to $r$ as the {\em row} of the node $(r,s)$ and $s$ as the {\em column} of the node $(r,s)$. 
Given nodes $A_1,\dots,A_k$ in $\la$, we denote 
$$\la_{A_1,\dots,A_k}:=\la\setminus\{A_1,\dots,A_k\}.$$
Given a node $A=(r,s)$ in row $r$ and column $s$, its {\em residue} is defined as 
$$
\res A:=s-r\pmod{p}\in I.
$$

\iffalse{
The {\em residue content} of a partition $\la$ is the tuple 
$\cont(\la):=(a_i)_{i\in I}$
such that $\la$ has exactly $a_i$ nodes of residue $i$ for each $i\in I$. For $j\in I$, let $\al_j$ be the tuple $(a_i)_{i\in I}$ with $a_i=\de_{i,j}$. We consider the tuples $(a_i)_{i\in I}$ as elements of $\Theta:=\sum_{i\in I}\Z\cdot\al_i$, the free $\Z$-module with basis $\{\al_i\mid i\in I\}$. 
Let 
\begin{equation*}\label{ETheta}
\Theta_n:=\bigg\{\theta=\sum_{i\in I}a_i\al_i\in\Theta\ \mid\ a_i\geq 0,\ \sum_{i\in I} a_i=n\bigg\}.
\end{equation*}
Given $\bi=(i_1,\dots,i_n)\in I^n$, we set 
$$
|\bi|:=\sum_{i\in I}\sharp\{r\mid i_r=i\}\cdot\al_i\in\Theta_n.
$$
%Partitions $\la,\mu\in\Par(n)$ have the same residue contents if and only if they have the same $p$-cores, see \cite[2.7.41]{JK}. 
}\fi

Let $i \in I$ and $\la\in\Par(n)$. 
A node $A 
\in \la$ (resp. $B\not\in\la$) of residue $i$ is called {\em $i$-removable} (resp. {\em $i$-addable}) for $\la$ 
if $\la_A$ (resp. $\la^B:=\la\cup\{B\}$) is a Young diagram of a partition. 
A node is called {\em removable} (resp. {\em addable}) if it is $i$-removable (resp. $i$-addable) for some $i$. 

Labeling the $i$-addable
nodes of $\la$ by $+$ and the $i$-removable nodes of $\la$ by $-$, the {\em $i$-signature} of 
$\la$ is the sequence of pluses and minuses obtained by going along the 
rim of the Young diagram from bottom left to top right and reading off
all the signs.
The {\em reduced $i$-signature} of $\la$ is obtained 
from the $i$-signature
by successively erasing all neighboring 
pairs of the form $-+$ (the result is independent of order in which one does this). 
%Note the reduced $i$-signature always looks like a sequence of $+$'s followed by $-$'s.
The nodes corresponding to  $-$'s in the reduced $i$-signature are
called {\em $i$-normal} for $\la$.
The leftmost $i$-normal  node is called {\em $i$-good} for $\la$. A node is called {\em normal} (resp. {\em good}) if it is $i$-normal (resp. $i$-good) for some $i$.
We denote 
$$
\eps_i(\la):=\sharp\{\text{$i$-normal nodes of $\la$}\}.
$$
There exists an $i$-good node for $\la$ if and only if $\eps_i(\la)>0$.  In this case we set 
$$
\tilde e_i \la:=\la_A
$$
where $A$ is the $i$-good node of $\la$.

Let $\la\in\Par_p(n)$. Note that the top removable node of $\la$ is always normal. If it is the only normal node of $\la$ we say that $\la$ is {\em Jantzen-Seitz}. The importance of Jantzen-Seitz partitions is explained by Lemma~\ref{Lemma39}(iv).

\subsection{Degenerate affine Hecke algebra}
\label{SSDAHA}
We denote by $\H_k$ the {\em degenerate affine Hecke algebra} given by generators $\s_1,\dots,\s_{k-1},\,\x_1,\dots,\x_k$ and relations 
\begin{align*}
\s_r^2=1,\ \ \s_r\s_{r+1}\s_r=\s_{r+1}\s_r\s_{r+1},\ \ \s_r\s_t=\s_t\s_r\ \text{for}\ |r-t|>1,
\\
\x_r\x_t=\x_t\x_r,\ \ \s_r\x_r=\x_{r+1}\s_r-1,\ \ \s_r\x_t=\x_t\s_r\ \text{for}\ t\neq r,r+1,
\end{align*}
see for example \cite[\S3.1]{KBook}. It is well-known that the subalgebra of $\H_k$ generated by the 
$\s_r$ is the group algebra $\F\Si_k$ and the subalgebra generated by the $\x_t$ is the polynomial algebra $\F[\x_1,\dots,\x_k]$, see for example \cite[Theorem 3.2.2]{KBook}. 

If $V\in\mod{\H_k}$, and $\ba=(a_1,\dots,a_k)\in \F^k$ is a $k$-tuple of scalars, we denote by $V_\ba$ the simultaneous generalized eigenspace for $\x_1,\dots,\x_k$ on $V$ corresponding to the eigenvalues $a_1,\dots,a_k$, respectively. We say that $V$ is {\em integral} if $V_\ba\neq 0$ implies $\ba\in I^n$ (recall that $I$ is identified with the prime subfield of $\F$). 

To every $k$-tuple $\ba=(a_1,\dots,a_k)\in \F^k$ one can associate an irreducible $\H_k$-module $L(\ba)=L(a_1,\dots,a_k)$ as in \cite[(5.14)]{KBook}. It follows easily from the definition that $L(\ba)_\ba\neq 0$. Moreover, by \cite[Lemma 4.2.2]{KBook}, $L(\ba)_\bb\neq 0$ only if the $k$-tuple $\bb$ is a permutation of the $k$-tuple $\ba$. 
Every irreducible $\H_k$-module is of this form but there might exist non-trivial isomorphisms $L(\ba)\cong L(\bb)$ if $\bb\neq\ba$ is a permutation of $\ba$. 

We now describe some of the irreducible modules $L(\ba)$ explicitly. 
For any $k$-tuple $\ba=(a_1,\dots,a_k)\in \F^k$, let $\De(\ba):=\Ind_{\F[\x_1,\dots,\x_k]}^{\H_k}\F_{\ba}$ where $\F_{\ba}$ is the $1$-dimensional $\F[\x_1,\dots,\x_k]$-module with each $\x_t$ acting with the scalar $a_t$. It is known that $\De(\ba)=L(\ba)$ in one of the following two cases: 
\begin{enumerate}
\item[{\rm (1)}] $a_1=a_2=\dots=a_k$, see \cite[Theorems 4.3.2]{KBook}
\item[{\rm (2)}] $|a_r-a_t|\neq 1$ for all $1\leq r\neq t\leq k$, see \cite[Theorems 6.1.4]{KBook}.
\end{enumerate}
We will need two more irreducible $\H_k$-modules:
$$L(a,a+1,\dots,a+k-1)\quad \text{and}\quad L(a,a-1,\dots,a-k+1).$$ Both are $1$-dimensional, with bases $\{v_\pm\}$, respectively, and the action is defined via
$$
\s_rv_\pm=\pm v_\pm\quad\text{and}\quad  
\x_tv_\pm=(a\pm(t-1))v_\pm.
$$

For the case $k=2$, we have constructed enough irreducible modules for a complete classification, see \cite[\S6.2]{KBook}. Namely, every irreducible $\H_2$-module is isomorphic to one of the modules $L(a,b)$ as above, and $L(a,b)\cong L(c,d)$ for $(c,d)\neq (a,b)$ if and only if $|a-b|\neq 0,1$ and $(c,d)=(b,a)$.

%the only non-trivial isomorphisms between these modules are $L(a,b)\cong L(b,a)$ if $|a-b|\neq 0,1$. 

We will also need the classification of integral irreducible $\H_3$-modules for the case $p=2$. Note that in this case $I=\{0,1\}$. In view of \cite[Lemma 6.2.1(ii)]{KBook}, 
$$
\{L(i,i,i),\ L(i,i,i+1),\ L(i,i+1,i),\,L(i+1,i,i)\mid i\in I\} 
$$
is then a complete and irredundant set of integral irreducible $\H_3$-modules up to isomorphism. We already know the modules $L(i,i,i)$ and $L(i,i+1,i)=L(i,i+1,i+2)$. As for $L(i,i,i+1)$ and $L(i+1,i,i)$, it will be sufficient to know that $L(i,i,i+1)=L(i,i,i+1)_{(i,i,i+1)}$ and $L(i+1,i,i)=L(i+1,i,i)_{(i+1,i,i)}$ are $2$-dimensional, with bases of the form $\{v,\s_1v\}$ and $\{w,\s_2w\}$, respectively. 

From the classification of irreducible modules in the previous two paragraphs we conclude:

\begin{Lemma} \label{LDeterm}
Suppose that either $k=2$, or $k=3$ and $p=2$. 
Let $i_1,\dots,i_k\in I$ and suppose that $L$ is an irreducible $H_k$-module with $L_{(i_1,\dots,i_k)}\neq 0$. Then $L\cong L(i_1,\dots,i_k)$. 
\end{Lemma}

The following lemma is immediate from the explicit constructions of irreducible modules above:

\begin{Lemma} \label{LHtoS}
We have for the restrictions of the irreducible $\H_k$-modules to the subalgebra $\F\Si_k$ generated by $\s_1,\dots,\s_{k-1}$: 
\begin{enumerate}
\item[{\rm (i)}] $\Res^{\H_k}_{\F\Si_k}L(a^k)\cong \F\Si_k$;
\item[{\rm (ii)}] if $|a_r-a_t|\neq 1$ for all $1\leq r\neq t\leq k$ then 
$\Res^{\H_k}_{\F\Si_k}L(a_1,\dots,a_k)\cong\F\Si_k$;
\item[{\rm (iii)}] $\Res^{\H_k}_{\F\Si_k}L(a,a+1,\dots,a+k-1)\cong\triv_{\Si_k}$; 
\item[{\rm (iv)}]$\Res^{\H_k}_{\F\Si_k}L(a,a-1,\dots,a-k+1)\cong\sgn_{\Si_k}$;
\item[{\rm (v)}] If $p=2$ then $\Res^{\H_3}_{\F\Si_3}L(i,i,i+1)\cong \Res^{\H_3}_{\F\Si_3}L(i+1,i,i)\cong D^{(2,1)}$.
\end{enumerate}

\end{Lemma}

\subsection{Modular branching rules}

The algebra $\H_k$ arises in the context of branching as follows. 
For $k=1,\dots,n$, let $H_{n,k}$ be the subalgebra of $\F\Si_n$ generated by the Murphy elements $L_{n-k+1},\dots,L_n$ and transpositions $(n-k+1,n-k+2),\dots,(n-1,n)$. It is well known that there is an algebra surjection 
$$\H_k\onto H_{n,k},\ \s_r\mapsto (n-k+r,n-k+r+1),\ \x_t\mapsto L_{n-k+t}.
$$
Moreover, it is well-known that the subalgebra $H_{n,k}\subseteq \F\Si_n$ centralizes the subalgebra $\F\Si_{n-k}\subseteq \F\Si_n$, and together they generate a subalgebra of $\F\Si_n$ isomorphic to $\F\Si_{n-k}\otimes H_{n,k}$. 
Now it is important that the restriction functor from $\F\Si_n$ to $\F\Si_{n-k}$ splits as the composition:
\begin{equation}\label{ERes}
\Res^{\F\Si_n}_{\F\Si_{n-k}}=
\Res^{\F\Si_{n-k}\otimes \F\Si_k}_{\F\Si_{n-k}}
\circ
\Res^{\F\Si_{n-k}\otimes H_{n,k}}_{\F\Si_{n-k}\otimes \F\Si_k}
\circ
\Res^{\F\Si_n}_{\F\Si_{n-k}\otimes H_{n,k}}.
\end{equation}

If the irreducible module $L(a_1,\dots,a_k)$ of $\H_k$ factors through its quotient $H_{n,k}$, we use the same notation $L(a_1,\dots,a_k)$ for the corresponding irreducible $H_{n,k}$-module. A necessary condition for that is that all $a_1,\dots,a_k\in I$, i.e. that $L(a_1,\dots,a_k)$ is integral, see \cite[Lemma 11.2.4]{KBook}.

For $V\in\mod{\F\Si_n}$ and $i\in I$ we denote by $e_iV$ the generalized $i$-eigenspace of $L_n$ on $V$. Since $L_n$ commutes with $\F\Si_{n-1}$, this is an $\F\Si_{n-1}$-submodule of $V$. It is well-known that $V\dar_{\Si_{n-1}}=\bigoplus_{i\in I}e_i V$, see for example \cite[Lemma 11.2.4]{KBook}. Hence
\begin{equation*}
V\dar_{\Si_{n-k}}=\bigoplus_{i_1\dots,i_k\in I}e_{i_1} \cdots e_{i_k}V.
\end{equation*}

\iffalse{
Suppose $V$ is in the block $B_\theta$ corresponding to $\theta\in \Theta_n$, for example if $V=D^\la$ with $\cont(\la)=\theta$. 
For $\eta \in\Theta_k$, denote by $\Res^\theta_{\theta-\eta,\eta}$ restricting to $\F\Si_{n-k}$ followed by projecting to the block $B_{\theta-\eta}$.
Then
\begin{equation}\label{EResEta}
\Res^\theta_{\theta-\eta,\eta}V=\bigoplus_{\bi=(i_1,\dots,i_k)\ \text{with}\ |\bi|=\eta}e_{i_1} \cdots e_{i_k}V
\end{equation}
and 
\begin{equation}\label{EResFirst}
\Res^{\F\Si_n}_{\F\Si_{n-k}\otimes H_{n,k}}=\bigoplus_{\eta \in\Theta_k}\Res^\theta_{\theta-\eta,\eta} V.
\end{equation}

The following lemma is now clear:

\begin{Lemma} \label{LResBlo}
Let $V$ be in the block corresponding to $\theta\in \Theta_n$. Then $\Res^\theta_{\theta-\eta,\eta}V$ is invariant under the action of  the subalgebra $\F\Si_{n-k}\otimes H_{n,k}\subseteq \F\Si_n$ and all composition factors of $\Res^\theta_{\theta-\eta,\eta} V$ are of the form $D^\mu\boxtimes L(\bi)$ with $\cont(\mu)=\theta-\eta$ and $|\bi|=\eta$. 
\end{Lemma}

%For example, if $\eta=\al_i+\al_j$ for $|i-j|\neq 0,1$, 
}\fi

\begin{Lemma} \label{LDeterm2}%{\rm \cite{}}%{\bf ()}
Suppose that either $k=2$, or $k=3$ and $p=2$. Let $i_1,\dots,i_k\in I$ and $\mu\in\Par_p(n-k)$. 
If $D^\mu$ is a composition factor of $e_{i_1}\cdots e_{i_k}V$ then $D^\mu\boxtimes L(i_1,\dots,i_k)$ is a composition factor of $\Res^{\F\Si_n}_{\F\Si_{n-k}\otimes H_{n,k}}V$. 
\end{Lemma}
\begin{proof}
Since $\Res^{\F\Si_n}_{\F\Si_{n-k}}=
\Res^{\F\Si_{n-k}\otimes  H_{n,k}}_{\F\Si_{n-k}}
\circ
\Res^{\F\Si_n}_{\F\Si_{n-k}\otimes H_{n,k}}$, 
if $D^\mu$ is a composition factor of $e_{i_1}\cdots e_{i_k}V$ then $D^\mu\boxtimes L$ is a composition factor of $\Res^{\F\Si_n}_{\F\Si_{n-k}\otimes H_{n,k}}$ for some irreducible $H_k$-module $L$ with $L_{(i_1,\dots,i_k)}\neq 0$. It remains to apply Lemma~\ref{LDeterm}. 
\end{proof}

We also define divided power functors $e_i^{(k)}$ as the space of the $\Si_k$-invariants on $e_i^kV$, which is isomorphic to taking $\Si_k$-anti-invariants on $e_i^kV$, see \cite[\S11.2,~\S8.3]{KBook} (in this situation we always use the embedding $\Si_k\into \Si_{n-k}\times\Si_k\into \Si_n$). 

We record some facts about the functors $e_i^{(k)}$ applied to irreducible modules for future reference. The following results are contained in \cite[Theorem 11.2.10]{KBook}, \cite[Theorem 1.4]{KDec}, see also \cite{KBrII}. Note that part (iv) follows from parts (ii) and (iii) but was first proved in \cite{JS} and \cite{k2}. 

\begin{Lemma}\label{Lemma39}
Let $V\in \mod{\F\Si_n}$, $\lambda\in\Parp(n)$, $i\in I$ and $k\in\Z_{\geq 0}$. Then:
\begin{enumerate}
\item[{\rm (i)}] $e_i^kV\cong(e_i^{(k)}V)^{\oplus k!}$;
\item[{\rm (ii)}]  $e_i^{(k)}D^\lambda\not=0$ if and only if $k\leq \eps_i(\lambda)$, in which case $e_i^{(k)}D^\lambda$ is a self-dual indecomposable module with socle and head both isomorphic to $D^{\tilde e_i^k\la}$.  
%\item[{\rm (iii)}]  $[e_i^{(k)}D^\lambda:D^{\tilde e_i^k\la}]=\binom{\eps_i(\lambda)}{k}$;
\item[{\rm (iii)}] Let $A$ be a removable node of $\la$ such that $\la_A$ is $p$-regular. Then $D^{\la_A}$ is a composition factor of $e_i D^\la$ if and only if $A$ is $i$-normal, in which case $[e_i D^\la:D^{\la_A}]$ is one more than the number of $i$-normal nodes for $\la$ above $A$. 
\item[{\rm (iv)}] $D^\la\dar_{\Si_{n-1}}$ is irreducible if and only if $\la$ is Jantzen-Seitz.
\end{enumerate}
\end{Lemma}

\subsection{The submodules $V_\pm$ and $V_\eps$}
\label{SSVPM}
Suppose first that $p\neq 2$. For a $V\in\mod{\F \Si_n}$, we denote by $V_\pm$ the $(\pm 1)$-eigenspace of $(n-1,n)$ on $V$. Note that $V_\pm\subseteq V$ is a submodule of $V\dar_{\Si_{n-2}}$ and 
$$
V\dar_{\Si_{n-2}\times\Si_2}=V_+\boxtimes \triv_{\Si_2}\,\oplus\,V_-\boxtimes \sgn_{\Si_2}.
$$
%From Lemma~\ref{LHtoS}, we deduce

\begin{Lemma} \label{LBr}
Let $p\neq2$, $V\in\mod{\F\Si_n}$, $\mu\in\Parp(n-2)$, and $i,j\in I$.  Suppose that $[e_ie_jV:D^\mu]\neq 0$.
If $j\neq i\pm1$ then $[V_\mp:D^\mu]\neq 0$.
\end{Lemma}
\begin{proof}
%We may assume that $V$ is in a fixed block corresponding to $\theta\in\Theta_n$. 
By Lemma~\ref{LDeterm2}, $D^\mu\boxtimes L(i,j)$ is a composition factor of $\Res^{\F\Si_n}_{\F\Si_{n-2}\otimes H_{n,2}}V$. If $j=i\pm1$,  acting with $(1\pm(n-1,n))$ on $\Res^{\F\Si_n}_{\F\Si_{n-2}\otimes H_{n,2}}V$ and taking into account 
Lemma~\ref{LHtoS}(iii)(iv), we deduce that $D^\mu$ is a composition factor of $(1\pm(n-1,n))V=V_\pm.$
If $j\neq i\pm1$, acting again with $(1\pm(n-1,n))$ on $\Res^{\F\Si_n}_{\F\Si_{n-2}\otimes H_{n,2}}V$ and taking into account 
Lemma~\ref{LHtoS}(i)(ii), we deduce that $D^\mu$ is a composition factor of $(1\pm(n-1,n))V=V_\pm.$
\end{proof}

%\subsection{The submodules $V_\eps$ for $p=2$}\label{SSp=2VZeta}
Now, let $p=2$. Denote by $\zeta$ a primitive $3$rd root of unity in $\F$. Then $\zeta^\eps\in \F$ makes sense for any $\eps\in\Z/3\Z=\{\0,\1,\2\}$. 
For $V\in\mod{\F \Si_n}$ and $\eps\in\Z/3\Z$, we denote by $V_\eps$ the $\zeta^\eps$-eigenspace of the $3$-cycle $(n-2,n-1,n)$ on $V$. Let $C_3$ be the cyclic subgroup generated by $(n-2,n-1,n)$, and denote by $\F(\eps)$ the $1$-dimensional module over $C_3$ with the generator acting by multiplying with $\zeta^\eps$. Then $V_\eps\subseteq V$ is a submodule of $V\dar_{\Si_{n-3}}$ and 
$$
V\dar_{\Si_{n-3}\times C_3}=V_{\0}\boxtimes \F(\0)\,\oplus\,V_{\1}\boxtimes \F(\1)\,\oplus\,V_{\2}\boxtimes \F(\2).
$$
Moreover, it is easy to see that $V_\1\cong V_\2$ as $\F\Si_{n-3}$-modules, the isomorphism coming from multiplying with $(n-1,n)$. 

\begin{Lemma} \label{LBr2}
Let $p=2$, $V\in\mod{\F\Si_n}$, $\mu\in\Parp(n-3)$, $i,j,k\in I$.  
Suppose that $[e_ie_je_kV:D^\mu]\neq 0$. If 
 $i=j$ or $j=k$ then $[V_\1:D^\mu]=[V_\2:D^\mu]\neq 0$. 
\end{Lemma}
\begin{proof}
By Lemma~\ref{LDeterm2}, $D^\mu\boxtimes L(i,j,k)$ is a composition factor of $\Res^{\F\Si_n}_{\F\Si_{n-3}\otimes H_{n,3}}V$. For $\eps=\1,\2$, acting with $(1+\zeta^\eps+\zeta^{2\eps})$ on $\Res^{\F\Si_n}_{\F\Si_{n-2}\otimes H_{n,2}}V$ and taking into account 
Lemma~\ref{LHtoS}(v), we deduce that $D^\mu$ is a composition factor of $(1+\zeta^\eps+\zeta^{2\eps})V=V_\eps.$
\end{proof}

\section{Rank and level}

\subsection{Definition and easy properties of rank} 
Let $V\in\mod{\F \Si_n}$. 

If $p\neq 2$, we define the {\em $2$-rank} of $V$ as the largest integer $r$ such that there exist distinct $i_1,j_1,\dots,i_r,j_r\in\{1,\dots,n\}$ and a simultaneous (non-zero) $(-1)$-eigenvector in $V$ for the transpositions $(i_1,j_1),\dots,(i_r,j_r)$.  
The $2$-rank of $V$ will be denoted $\rank_2(V)$. 
Equivalently, consider the 
Young subgroup
\begin{equation}\label{EE2}
E_2 := \Si_2 \times \ldots \times \Si_2 \times \Si_{n-2m}\leq \Si_n,
\end{equation} 
where $m := \lfloor n/2 \rfloor$ (so $n-2m=0$ or $1$, and $ \Si_{n-2m}$ is always the trivial group). 
The irreducible $\F E_2$-modules are of the form $L=L_1\boxtimes\dots\boxtimes L_m\boxtimes L_{m+1}$ (where $L_{m+1}=\triv_{\Si_{n-2m}}$). Define the {\em type}\, of such an irreducible $\F E_2$-module $L$ as 
$$\type_{E_2}(L):=\sharp\{r\mid 1\leq r\leq m\ \text{and}\ L_r\cong\sgn_{\Si_2}\}.$$ 
The {\em type} $\type_{E_2}(W)$ of any $W\in\mod{\F E_2}$ is defined to be the maximal {\em type} of its composition factors. Now, for $V\in\mod{\F \Si_n}$, we have that 
\begin{equation}\label{ERT2}
\rank_2(V)=\type_{E_2}(V\dar_{E_2}). 
\end{equation}

If $p\neq 3$, let $\zeta$ be a primitive $3$rd root of unity in $\F$ (here we assume that $\F$ does contain a primitive $3$rd root of unity), and define the {\em $3$-rank} of $V$ as the largest integer $r$ such that there exist distinct $i_1,j_1,k_1,\dots,i_r,j_r,k_r\in\{1,\dots,n\}$ and a simultaneous (non-zero) $\zeta$-eigenvector in $V$ for the $3$-cycles $(i_1,j_1,k_1),\dots,(i_r,j_r,k_r)$. The $3$-rank of $V$ will be denoted $\rank_3(V)$. Equivalently, consider the 
Young subgroup
\begin{equation}\label{EE3}
E_3 := \Si_3 \times \ldots \times \Si_3 \times \Si_{n-3l}\leq \Si_n,
\end{equation} 
where $l := \lfloor n/3 \rfloor$ (so $n-3l\leq 2$). 
The irreducible $\F E_3$-modules are of the form $L=L_1\boxtimes\dots\boxtimes L_l\boxtimes L_{l+1}$. 
Define the {\em type}\, of such an irreducible $\F E_3$-module $L$ as 
$$\type_{E_3}(L):=\sharp\{r\mid 1\leq r\leq l\ \text{and}\ L_r\cong D^{(2,1)}\}.$$ 
The {\em type} $\type_{E_3}(W)$ of any $W\in\mod{\F E_3}$ is defined to be the maximal {\em type} of its composition factors. Now, for $V\in\mod{\F \Si_n}$, we have that 
\begin{equation}\label{ERT3}
\rank_3(V)=\type_{E_3}(V\dar_{E_3}). 
\end{equation}
Note that this description of $3$-rank does not require the assumption that $\F$ contains a primitive $3$rd root of unity.

%Define the rank of such an irreducible $\F E_3$-module is $\sharp\{r\mid 1\leq r\leq l\ \text{and}\ L_r\cong D^{(2,1)}\}$. The rank of any $W\in\mod{\F E_3}$ is defined to be the maximal rank of its composition factors. Now, for $V\in\mod{\F \Si_n}$, we have that $\rank_3(V)$ equals the rank of the $\F E_3$-module $V\dar_{E_3}$. 

\iffalse{
For an $\F\Si_n$-module $V$, we define its {\em rank} to be 
$$
\rank(V):=
\left\{
\begin{array}{ll}
\rank_2(V) &\hbox{if $p\neq 2$,}\\
\rank_3(V) &\hbox{if $p=2$.}
\end{array}
\right.
$$
}\fi

The following properties are immediate:

\begin{Lemma} \label{LRed}
Let $V$ be a reduction modulo $p$ of a $\C \Si_n$-module $V_\C$. Then $\rank_2(V)=\rank_2(V_\C)$ if $p\neq 2$, and $\rank_3(V)=\rank_3(V_\C)$ if $p\neq 3$.  
\end{Lemma}

\begin{Lemma} \label{L1}
If $0=W_0\subseteq W_1\subseteq\dots\subseteq W_k= V$ is a chain of $\F \Si_n$-submodules then 
\begin{align*}
\rank_2(V)&=\max\{\rank_2(W_i/W_{i-1})\mid 1\leq i\leq k\}\qquad\text{if $p\neq 2$},\\
\rank_3(V)&=\max\{\rank_3(W_i/W_{i-1})\mid 1\leq i\leq k\}
\qquad\text{if $p\neq 3$}. 
\end{align*}
\end{Lemma}

Let $V\in\mod{\F \Si_n}$. 
If $p\neq 2$, since $r$-tuples of non-overlapping transpositions are all conjugate in $\Si_n$, we have that $\rank_2(V)$ is the largest integer $r$ such that the transpositions 
$(n,n-1),(n-2,n-3),\dots, (n-2r+2,n-2r+1)$ have a simultaneous $(-1)$-eigenvector in $V$. Similarly, if $p\neq 3$, we have that $\rank_3(V)$ is the largest integer $r$ such that 
the $3$-cycles  
$(n,n-1,n-2),(n-3,n-4,n-5),\dots, (n-3r+3,n-3r+2,n-3r+1)$ have a simultaneous $\zeta$-eigenvector in $V$. 
Recall the notation $V_-$ and $V_\1$ from \S\ref{SSVPM}. 
% and \S\ref{SSp=2VZeta}. 
The description of $\rank_2(V)$ and $\rank_3(V)$ given in this paragraph makes the  following clear:

\begin{Lemma} \label{L2}
Let $V$ be an $\F \Si_n$-module. If $p\neq 2$ then $\rank_2(V)=1+\rank_2(V_-)$. If $p\neq 3$ then $\rank_3(V)=1+\rank_3(V_\1)$.
\end{Lemma}

\iffalse{
\begin{Lemma} %\label{}%{\rm \cite{}}%{\bf ()}
Let $V$ be an $\F \Si_n$-module. If $p\neq 2$ then $\max\big(\rank_2(V),\rank_2(V\otimes \sgn)\big)=\lfloor n/2\rfloor$. 
\end{Lemma}
\begin{proof}
By Lemma~\ref{L1}, we may assume that $V$ is irreducible and different from $\triv$ and $\sgn$, so that $V_\pm\neq 0$. We apply induction on $n$, the cases $n=1,2$ being trivial. Let $n>2$. 
We may assume that $r_2(V)<\lfloor n/2\rfloor$. Then by Lemma~\ref{L2} we have $\rank_2(V_-)<\lfloor n/2\rfloor-1=\lfloor(n-2)/2\rfloor$. 
So by the inductive assumption we have $\rank_2(V_-\otimes\sgn_{\Si_{n-2}})=\lfloor(n-2)/2\rfloor=\lfloor n/2\rfloor-1$.  Note that $V_-\otimes \sgn_{\Si_{n-2}}\cong 
(V\otimes\sgn_{\Si_n})_+$ as $\F\Si_{n-2}$-modules. 
 So 
\end{proof}
}\fi

For an irreducible $\F\Si_n$-module $D$, we define its {\em rank} to be 
$$
\rank(D):=
\left\{
\begin{array}{ll}
\min\big(\rank_2(D),\rank_2(D\otimes\sgn)\big) &\hbox{if $p\neq 2$,}\\
\rank_3(D) &\hbox{if $p=2$.}
\end{array}
\right.
$$
For an arbitrary $V\in\mod{\F\Si_n}$, we define its rank to be 
$$
\rank(V):=\max\{\rank(D)\mid D\ \text{is a composition factor of}\ V\}.
$$
It is clear that
\begin{equation}\label{ERankSign}
\rank(V)=\rank(V\otimes\sgn).
\end{equation}

\begin{Remark} %\label{}%{\rm \cite{}}%{\bf ()}
{\rm 
Somewhat surprisingly, in characteristic $0$, it follows easily from Proposition~\ref{PSpechtRank}(i) below that $\max\big(\rank_2(V),\rank_2(V\otimes\sgn)\big)=\lfloor n/2\rfloor$ for all  $V\in\mod{\F\Si_n}$. This is not true in positive characteristic in general. For example, by the description of the Mullineux bijection in \cite{FK}, we have $D^\la\otimes \sgn\cong D^\la$ for $\la=(7,3,2)$ in characteristic $p=3$. On the other hand, by Theorem~\ref{TRL} below, we have that $\rank_2(D^\la)=5$ in this case. 
\iffalse{
We make some additional observation about the module  $V := D^{(7,3,2)}$ in characteristic $3$. Note that by definition of $\rank_2$, the composition factors of $V\dar_{E_2}$ must have ranks $0\leq a \leq 5$, and the 
composition factors of $(V \otimes \sgn)\dar_{E_2}$ have ranks $6-a$ for those $a$ that do occur as ranks of composition factors of $V\dar_{E_2}$. But since 
$V \otimes \sgn \cong V$ has $2$-rank equal to $5$, this shows that $V\dar_{E_2}$ does not have a composition factor of rank $0$. Thus, in the modular case, one cannot hope for ``continuity'' of $\rank_2$, that is, we cannot hope that $V\dar_{E_2}$ will always  afford composition factors of all the 
intermediate ranks between $0$ and $\rank_2(V)$. 
\color{red} Does continuity hold in characteristic $0$? I forgot the answer to that. If it does, it needs to be mentioned here, perhaps with a quick explanation why. Is it true by Michael's induction argument? If yes, it might worth making it a lemma somewhere.\color{black}
}\fi
Taking $W = \sgn$, this example also shows that, in general, $\rank_2(V \otimes W)$ is {\it not} equal to
$\min\bigl( \rank_2(V)+\rank_2(W),\lfloor n/2 \rfloor \bigr)$, cf. Theorem~\ref{TTensRank}(i) below. (A more boring example is given by considering $\sgn\otimes\sgn$.)
}
\end{Remark}

\subsection{Some applications of rank} 

In this subsection we establish some results which explain why rank is of interest for us. First, we have lower bounds for dimensions in terms of rank:

\begin{Theorem} \label{TDimRank}
Let $V$ be an $\F\Si_n$-module and $r=\rank(V)$. 
\begin{enumerate}
\item[{\rm (i)}] If $p\neq 2$ then $\dim V\geq {\lfloor n/2\rfloor\choose r}$. 
\item[{\rm (ii)}] If $p= 2$ then $\dim V\geq 2^r{\lfloor n/3\rfloor\choose r}$. 
\end{enumerate}
\end{Theorem}

\begin{proof}
(i) By (\ref{ERT2}), $V\dar_{E_2}$ or $(V\otimes\sgn)\dar_{E_2}$ contains a composition factor $L$ with $\type_{E_2}(L) =r$. 
Let $N_2:=N_{\Si_n}(E_2) \cong E_2 \rtimes \Si_{m}$ where $m = \lfloor n/2 \rfloor$. Then $L$ is a direct summand of $M\dar_{E_2}$, where $M$ is a composition factor of $V\dar_{N_2}$ or $(V\otimes\sgn)\dar_{N_2}$. Considering the $N_2$-orbit on the irreducible summands of $M\dar_{E_2}$, it follows that $M$ has dimension at least $\binom{m}{r}$ (in fact, exactly $\binom{m}{r}$ since in this case we must have $M\cong\Ind_{E_2}^{N_2}L$). 

(ii) By (\ref{ERT3}), $V\dar_{E_3}$ contains a composition factor $L$ with $\type_{E_3}(L) =r$. Note that $\dim L=2^r$. 
Let $N_3:=N_{\Si_n}(E_3) \cong E_3 \rtimes \Si_{l}$ where $l = \lfloor n/3 \rfloor$. Then $L$ is a direct summand of $M\dar_{E_3}$, where $M$ is a composition factor of $V\dar_{N_3}$. Considering the $N_3$-orbit on the irreducible summands of $M\dar_{E_3}$, it follows that $M$ has dimension at least $\binom{l}{r}\dim L$. 
\end{proof}

We also have the following `additivity' property of ranks with respect to tensor products:  

\begin{Theorem} \label{TTensRank}
Let $V,W$ be $\F\Si_n$-modules. 
\begin{enumerate}[\rm(i)]
\item If $p \neq 2$ and $\rank_2(V)+\rank_2(W)\leq n/2$ then  $\rank_2(V\otimes W)=\rank_2(V)+\rank_2(W)$.  
\item If $p \neq 2$ and $\rank(V)+\rank(W)\leq n/4$ then  $\rank(V\otimes W)=\rank(V)+\rank(W)$. 
\item If $p = 2$ and $\rank(V)+\rank(W)\leq n/3$ then  $\rank(V\otimes W)=\rank(V)+\rank(W)$.  
\end{enumerate}
\end{Theorem}

\begin{proof}
Suppose first that $p \neq 2$.
By definition of the rank, we may assume that $V$ and $W$ are irreducible. Moreover, in view of (\ref{ERankSign}), tensoring with $\sgn$ if necessary, we may assume that $\rank(V)=\rank_2(V) =: r$ and $\rank(W)=\rank_2(W) =: s$. 

(i) By assumption,
$r+s \leq m:= \lfloor n/2 \rfloor$. By definition of $\rank_2$, we know that there is some $0 \neq v \in V$ such that $(2i-1,2i)$ sends $v$ to $-v$ when $1 \leq i \leq r$ and
fixes $v$ if $r+1 \leq i \leq m$. Similarly, there is some $0 \neq w \in W$ such that $(2i-1,2i)$ fixes $w$ when $1 \leq i \leq m-s$ and  
sends $v$ to $-v$ if $m-s+1 \leq i \leq m$. As $r+s \leq m$, 
the vector $v\otimes w\in V\otimes W$ is a $(-1)$-eigenvector for the $r+s$ non-overlapping transpositions $\{(2i-1,2i)\mid 1 \leq i \leq r\ \text{or}\ m-s+1 \leq j \leq m\}$, 
showing 
$\rank_2(V \otimes W) \geq r+s.$

To bound $\rank_2(V \otimes W)$ from the above, note that the composition factors of $(V\otimes W)\dar_{E_2}$ are of the form 
\begin{equation}\label{ECF}
M(L,L'):=(L_1\otimes L_1')\boxtimes\dots\boxtimes (L_m\boxtimes L_m')\boxtimes (L_{m+1}\otimes L_{m+1}'),
\end{equation}
where $L=L_1\boxtimes\dots\boxtimes L_m\boxtimes L_{m+1}$ is a composition factor of
$V\dar_{E_2}$, and $L'=L_1'\boxtimes\dots\boxtimes L_m'\boxtimes L_{m+1}'$ is a composition factor of $W\dar_{E_2}$. 
Let $A_L=\{i\mid 1\leq i\leq m\ \text{and}\ L_i\cong\sgn\}$ and $A_{L'}=\{i\mid 1\leq i\leq m\ \text{and}\ L_i'\cong\sgn\}$. Note that $|A_L|\leq r$ and $|A_{L'}|\leq r'$. 
Now, for any $1\leq i\leq m$, we have $L_i \otimes L_i'\cong \sgn$  if and only if 
$i \in (A_L \smallsetminus A_{L'}) \cup (A_{L'} \smallsetminus A_L)$. Thus the composition factor $M(L,L')$ as in (\ref{ECF}) has  $\type_{E_2}(M(L,L'))\leq |A_L|+|A_{L'}| \leq r+r'$,
and hence by (\ref{ERT2}), we have 
$\rank_2(V \otimes W) \leq r+r',$ proving (i). 

(ii) By assumption, $r+s \leq n/4$, i.e. $r+s \leq k:= \lfloor n/4 \rfloor$. Then $n \geq 4k$ and $m \geq 2k$. 
The above analysis shows that 
$$\type_{E_2}(M(L,L')\otimes\sgn)=m-|A_L \smallsetminus A_{L'}|- |A_{L'} \smallsetminus A_L| \geq m-(r+s) \geq k \geq r+s.$$
So by (\ref{ERT2}), for any composition factor $D$ of $V \otimes W \otimes \sgn$ we have $\rank_2(D) \geq r+s$. On the other hand, by (i), 
for a composition factor $D \otimes \sgn$ of $V \otimes W$, we have $\rank_2(D \otimes \sgn) \leq r+s$. Hence
$\rank(D \otimes \sgn) = \rank_2(D \otimes \sgn)$. It follows that $\rank(V \otimes W) = \rank_2(V \otimes W)$, and the latter equals
$r+s$ by (i). 

(iii) The proof is similar to the proof of (i) using that, for irreducible $F\Si_3$-modules $L_i,L_i'$, we have that $D^{(2,1)}$ appears as a composition factor of $L_i\otimes L_i'$ unless $\dim L_i=\dim L_i'=1$. 
\end{proof}

\iffalse{
\color{red}
{\bf Note from Sasha:}

Can we slightly improve on Theorem~\ref{TTensRank} to get rather the following statement:

(*) Let $V,W$ be $\F\Si_n$-modules. Then $\rank(V\otimes W)=\min(\rank(V)+\rank(W),\lfloor n/2\rfloor)$.  

\color{black}
}\fi

\subsection{Level of a partition and rank of a Specht module}
\label{SSLRS}

For a partition $\la=(\la_1,\la_2,\dots)$ of $n$ we define the {\em level} of $\la$ to be 
$$
\level(\la):=n-\la_1.
$$
%and call it the {\em level} of $\la$. 

We can compute ($k$-)ranks of Specht modules in terms of levels as follows: 

\begin{Proposition}\label{PSpechtRank}%{\rm \cite{}}%{\bf ()}
Let $\la\in\Par(n)$. 
\begin{enumerate}
\item[{\rm (i)}] If $p=0$ then $\rank_2(S^\la)=\min\{\level(\la),\lfloor n/2\rfloor\}$.
\item[{\rm (ii)}] If $p=0$ then $\rank_3(S^\la)=\min\{\level(\la),\level(\la'),\lfloor n/3\rfloor\}.$
\item[{\rm (iii)}] If $p\neq 2$ then $\rank_2(S^\la)=\min\{\level(\la),\lfloor n/2\rfloor\}.$
\item[{\rm (iv)}] If $p\neq 3$ then $\rank_3(S^\la)=\min\{\level(\la),\level(\la'),\lfloor n/3\rfloor\}.$
\end{enumerate}
 
\end{Proposition}
\begin{proof}
In view of Lemma~\ref{LRed}, (iii) follows form (i) and (iv) follows from (ii). Thus it suffices to prove (i) and (ii). 

(i) We apply induction on $n$, the induction base case $n=1$ being clear. If $n=2$ the result is also clear since $S^{(2)}$ is the trivial module and $S^{(1,1)}$ is the sign module.  Let $n>2$. We may assume that $\la\neq (n)$ since in this case we have $S^{(n)}\cong\triv$, and the result clearly holds: $\level((n))=0=\rank_2(\triv)$

Denote by $P$ the set of partitions $\mu$ of $n-2$ such that the Young diagram of $\mu$ is obtained from that of $\la$ by removing two nodes {\em not}\, in the same row. 
By (a special case of) Littlewood-Richardson rule \cite[2.8.13]{JK}, for $\mu\in\Par(n-2)$, the module $S^\mu$ appears as a composition factor of $S^\la_-$ if and only if $\mu\in P$. 

By Lemmas~\ref{L1},\ref{L2}, the previous paragraph, and inductive assumption, it suffices to prove that 
\begin{equation}\label{EMaxMin}
\min\{\level(\la),\lfloor n/2\rfloor\}=1+\max\big\{\min\{\level(\mu),\lfloor (n-2)/2\rfloor\}\mid \mu\in P\big\}.
\end{equation}

%In proving (\ref{EMaxMin}) we will use the following observation. Suppose $\mu\in P$ is obtained from $\la$ by removing nodes $A$ and $B$ (not in the same row). Then:

Suppose first that $\level(\la)\leq\lfloor n/2\rfloor$, so that the left hand side above equals $\level(\la)$. In the special case where $n$ is even and $\la$ is the two row partition $(n/2,n/2)$, we have that $P=\{(n/2-1,n/2-1)\}$ and the equality (\ref{EMaxMin}) is clear. In all other cases, the node $A:=(1,\la_1)$ in the first row is removable. If there are no other removable nodes in $\la$ then $\la=(n)$---the case we have already excluded. So we may assume that $\la$ has a removable node $B\neq A$. So $\mu:=\la_{A,B}\in P$ and $\level(\mu)=\level(\la)-1\leq\lfloor n/2\rfloor-1=\lfloor (n-2)/2\rfloor$. So $\min\{\level(\mu),\lfloor (n-2)/2\rfloor\}=\level(\la)-1$. On the other hand, any $\nu\in P$ is obtained either by removing $A$ and some node $C$ from a row $>1$, or two nodes $C$ and $D$ from rows $>1$. In the first case we have $\min\{\level(\nu),\lfloor (n-2)/2\rfloor\}=\level(\la)-1$ and in the second case we have $\min\{\level(\nu),\lfloor (n-2)/2\rfloor\}=\level(\la)-2$. This completes the proof of (\ref{EMaxMin}) in the case $\level(\la)\leq\lfloor n/2\rfloor$. 

On the other hand, if $\level(\la)>\lfloor n/2\rfloor$ then we have $\min\{\level(\mu),\lfloor (n-2)/2\rfloor\}=\lfloor n/2\rfloor-1$ for all $\mu\in P$,  which implies (\ref{EMaxMin}) in this case. 

(ii) Note that $S^{\la'}\cong S^\la\otimes\sgn$ by \cite[Theorem 6.7]{JamesBook}, so $\rank_3(S^{\la'})=\rank_3(S^\la)$, and we may assume that $\level(\la)\leq \level(\la')$. We apply induction on $n$, the induction base case $n=1$ being clear. If $n=2$ or $3$ the result is also easy to check.  Let $n>3$. We may also assume that $\la\neq (n)$ since in this case we clearly have $\rank(S^{(n)})=\level(\la)=0$. 

Denote now by $P$ the set of partitions $\mu$ of $n-3$ such that the Young diagram of $\mu$ is obtained from that of $\la$ by removing three nodes {\em neither}\, all in the same row {\em nor} \, all in the same column. 
By (a special case of) Littlewood-Richardson rule \cite[2.8.13]{JK}, for $\mu\in\Par(n-3)$, the module $S^\mu$ appears as a composition factor of $S^\la_\1$ if and only if $\mu\in P$. 

By Lemmas~\ref{L1}, \ref{L2}, the previous paragraph, and inductive assumption, it suffices to prove that 
\begin{equation}\label{EMaxMin2}
\min\{\level(\la),\lfloor n/3\rfloor\}=1+\max\big\{\min\{\level(\mu),\level(\mu'),\lfloor (n-3)/3\rfloor\}\mid \mu\in P\big\}.
\end{equation}

%In proving (\ref{EMaxMin}) we will use the following observation. Suppose $\mu\in P$ is obtained from $\la$ by removing nodes $A$ and $B$ (not in the same row). Then:

Suppose first that $\level(\la)\leq\lfloor n/3\rfloor$, so that the left hand side above equals $\level(\la)$. Note that the assumption $\level(\la)\leq\lfloor n/3\rfloor$ guarantees that $A:=(1,\la_1)$ is removable for $\la$ and $B:=(1,\la_1-1)$ is removable for $\la_A$. Since we have assumed that $\la\neq (n)$, the partition $\la$ must have a removable node $C$ not in row $1$. So $\mu:=\la_{A,B,C}\in P$ and $\level(\mu)=\level(\la)-1\leq\lfloor n/3\rfloor-1=\lfloor (n-3)/3\rfloor$. The assumption $\level(\la)\leq\lfloor n/3\rfloor$ easily implies that $\level(\mu)\leq \level(\mu')$, so 
 $\min\{\level(\mu),\level(\mu'),\lfloor (n-3)/3\rfloor\}=\level(\la)-1$. 
 On the other hand, it is clear that $\level(\nu)\leq \level(\la)-1$ for any $\nu\in P$. This completes the proof of (\ref{EMaxMin}) in the case $\level(\la)\leq\lfloor n/3\rfloor$. 

Suppose that $\level(\la)>\lfloor n/3\rfloor$. Let $\mu\in P$. Then $\min\{\level(\mu),\level(\mu')\}\geq \level(\la)-3$, so taking into account that $\lfloor (n-3)/3\rfloor=\lfloor n/3\rfloor-1$, we have  $\min\{\level(\mu),\level(\mu'),\lfloor (n-3)/3\rfloor\}=\lfloor n/3\rfloor-1$ unless 
$\min\{\level(\mu),\level(\mu')\}= \level(\la)-3$ and $\level(\la)=\lfloor n/3\rfloor+1$. 
But it is easy to see that in the case $\level(\la)=\lfloor n/3\rfloor+1$ there always exists $\mu\in P$ with $\min\{\level(\mu),\level(\mu')\}\geq \level(\la)-2$. This completes the proof of the case $\level(\la)>\lfloor n/3\rfloor$. 
\end{proof}

\subsection{Ranks of irreducible modules }\label{SLevel}
In this subsection we relate ranks and levels for irreducible representations of symmetric groups in positive characteristic, the characteristic zero case being covered by Proposition~\ref{PSpechtRank}(i),(ii). So throughout the subsection we may assume that $p>0$. Suppose $p\neq 2$. Then by definition 
$\rank(D^\la)=\min\left(\rank_2(D^\la),\rank_2(D^\la\otimes \sgn)\right)$, and $D^\la\otimes \sgn$ is in principle known, see \cite{FK}, so it suffices to compute $\rank_2(D^\la)$. If $p=2$ then by definition $\rank(D^\la)= \rank_3(D^\la)$ so we compute $\rank_3(D^\la)$ in this case.

\begin{Theorem} \label{TRL}
Let $p\neq 2$ and $\la\in\Parp(n)$. Then 
$\rank_2(D^\la)=\min\{\level(\la),\lfloor n/2\rfloor\}.$ 
\end{Theorem}
\begin{proof}
We may assume that $\la\neq (n)$ since in that case we have $D^{(n)}=\triv$ and the result is clear. 

We apply induction on $n$, the cases $n=1,2$ being clear. Let $n>2$ and denote by $P$ the set of $p$-regular partitions of $n-2$ such that $\{D^\mu\mid \mu\in P\}$ is the set of composition factors of $D^\la_-$. 

Since $D^\la$ is a composition factor of $S^\la$, Lemma~\ref{L1} and  Proposition~\ref{PSpechtRank}(iii) imply that  $\rank_2(D^\la)\leq \min\{\level(\la),\lfloor n/2\rfloor\}.$ So by Lemmas~\ref{L1},\,\ref{L2} and the inductive assumption it suffices to prove that there is $\mu\in P$ such that 
$$
\min\{\level(\mu),\lfloor (n-2)/2\rfloor\}=\min\{\level(\la),\lfloor n/2\rfloor\}-1.
$$

{\em Case 1:} $\level(\la)>\lfloor n/2\rfloor$. 

It suffices to find $\mu\in P$ with $\level(\mu)\geq \level(\la)-2$ since then 
$$\level(\mu)\geq \level(\la)-2\geq \lfloor n/2\rfloor-1=
\lfloor (n-2)/2\rfloor\}$$ 
and 
$$\min\{\level(\mu),\lfloor (n-2)/2\rfloor\}=\lfloor (n-2)/2\rfloor=\lfloor n/2\rfloor-1=\min\{\level(\la),\lfloor n/2\rfloor\}-1,$$
as required. 

Let $A$ be the top removable node of $\la$ such that $\la_A$ is $p$-regular, and let $i=\res A$. 

If $A$ has normal nodes above it then they are all necessarily of residue $i$, so $\eps_i(\la)\geq 2$ and by 
Lemma~\ref{Lemma39}(iii), we have that $D^{\tilde e_i^2\la}$ is a composition factor of 
$e_i^2D^\la$. By Lemma~\ref{LBr}, we can now take $\mu=\tilde e_i^2\la$. 

If $A$ does not have normal nodes above it and $A$ is not in row $1$, then the node $B$ just above $A$ is normal for $\la_A$ and has residue $i+1$. So there exists an $(i+1)$-good node $C$ for $\la_A$. We deduce from Lemma~\ref{Lemma39}(iii) that $D^{\la_{A,C}}$ is a composition factor of $e_{i+1}e_iD^\la$. By Lemma~\ref{LBr}, we can now take $\mu=\la_{A,C}$. 

If $A$ is in row $1$ then, since $\la\neq (n)$, there exists a second removable node from the top, call it $B$. Let $j=\res B$. If $j\neq i+1$ then $B$ is normal for $\la$. Let $C$ be the $j$-good node for $\la$. Then $A$ is normal for $\la_C$ so by 
Lemma~\ref{Lemma39}(iii), we have that $D^{\la_{A,C}}$ is a composition factor of $e_ie_jD^\la$, and by Lemma~\ref{LBr}, we can take $\mu=\la_{A,C}$. If $j=i+1$ then $B$ is normal for $\la_A$. Let $E$ be the $(i+1)$-good node for $\la_A$. Then by 
Lemma~\ref{Lemma39}(iii),  we have that $D^{\la_{A,E}}$ is a composition factor of $e_{i+1}e_iD^\la$, and by Lemma~\ref{LBr}, we can take $\mu=\la_{A,E}$.

{\em Case 2:}  $\level(\la)\leq\lfloor n/2\rfloor$. 

We have to find 
$\mu\in P$ such that $\min\{\level(\mu),\lfloor (n-2)/2\rfloor\}=\level(\la)-1$. 
In the exceptional case 
where $n$ is even and $\la$ is the two row partition $(n/2,n/2)$, we have $\level(\la)=n/2$ and 
by Lemmas~\ref{Lemma39}(iii) and \ref{LBr}, 
we can take $\mu=(n/2-1,n/2-1)$. So we may assume that the node $A:=(1,\la_1)$ is removable for $\la$. Let $A$ have residue $i$.

Assume first that $\la$ is not Jantzen-Seitz. In this case there exists a good node $B$ of $\la$ different from $A$. Set $\mu:=\la_{A,B}$. Note that $\level(\mu)=\level(\la)-1$, so it suffices to prove that $\mu\in P$. Let $j$ be the residue of $B$. If $j=i$, then by Lemma~\ref{Lemma39}(iii), we have that $D^{\mu}$ is a composition factor of $e_i^2 D^\la$, so $\mu\in P$ by Lemma~\ref{LBr}. If $j\neq i$ then by Lemma~\ref{Lemma39}(iii),  we have that $D^\mu$ is a composition factor of $e_ie_jD^\la$ which implies $\mu\in P$ by Lemma~\ref{LBr} unless $j=i+1$. If $j=i+1$ then $B$ is a normal node for $\la_A$ so by Lemma~\ref{Lemma39}(iii), we have that $D^\mu$ is a composition factor of $e_{i+1}e_iD^\la$ which again implies $\mu\in P$ by Lemma~\ref{LBr}. 

Assume that $\la$ is Jantzen-Seitz. Let $B$ be the second removable node of $\la$ from the top. Note that $B$ has residue $i+1$. So setting $\mu:=\la_{A,B}$, by Lemma~\ref{Lemma39}(iii), we have that $D^\mu$ is a composition factor of $e_{i+1}e_i D^\la$, which implies $\mu\in P$ by Lemma~\ref{LBr}.
\end{proof}

\begin{Theorem} \label{TRL2}
Let $p= 2$ and $\la\in\Par_2(n)$. Then 
$\rank_3(D^\la)=\min\{\level(\la),\lfloor n/3\rfloor\}.$ 
\end{Theorem}
\begin{proof}
We may assume that $\la\neq (n)$ since in that case we have $D^{(n)}=\triv$ and the result is clear. So we may assume that $\level(\la)>0$. 

We apply induction on $n$ the cases $n=1,2,3$ being clear. Let $n>3$ and denote by $P$ the set of $p$-regular partitions of $n-3$ such that $\{D^\mu\mid \mu\in P\}$ is the set of composition factors of $D^\la_\1$. %To get information on $P$ we will repeatedly use Lemmas~\ref{LBr2} and \ref{Lemma39}(iii) without further comment. 

Consider the nodes $A:=(1,\la_1)$ and of residue $i:=\la_1-1\pmod{2}$ and $B:=(2,\la_2)$ of residue $j:=\la_2-2\pmod{2}$. Since $\la$ is $2$-regular, $A$ and $B$ are  removable for $\la$. Let us also consider the node $C:=(1,\la_1-1)$ has residue $i-1=i+1$. 
%Now, if $D^{\mu}$ is a composition factor of either $ 

Suppose that $j=i$. Then $\eps_i(\la)\geq 2$ and $C$ is an $i+1$-normal node for $\tilde e_i^2\la$. So by Lemma~\ref{Lemma39}(iii) we have that $D^{(\tilde e_i^2\la)_C}$ is a composition factor of $e_{i+1}e_i^2 D^\la$, so by Lemmas~\ref{LBr2} we deduce:
\begin{equation}\label{EThx1}
(\tilde e_i^2\la)_C\in P\qquad \text{if $i=j$}.
\end{equation}
On the other hand, if $j=i+1$ then $\la_A$ is $2$-regular and $\eps_{i+1}(\la_A)\geq 2$. So by Lemma~\ref{Lemma39}(iii) we have that $D^{\tilde e_{i+1}^2(\la_A)}$ is a composition factor of $e_{i+1}^2e_i D^\la$, and by Lemmas~\ref{LBr2} we deduce:
\begin{equation}\label{EThx2}
\tilde e_{i+1}^2(\la_A)\in P\qquad \text{if $i=j+1$}.
\end{equation}

Since $D^\la$ is a composition factor of $S^\la$, Lemma~\ref{L1} and  Proposition~\ref{PSpechtRank}(iv) imply  
$$
\rank_3(D^\la)\leq 
\min\{\level(\la),\level(\la'),\lfloor n/3\rfloor\}\leq 
\min\{\level(\la),\lfloor n/3\rfloor\}.
$$ 
So by Lemmas~\ref{L1},\,\ref{L2} and the inductive assumption it suffices to prove that there is $\mu\in P$ such that 
\begin{equation}\label{E}
\min\{\level(\mu),\lfloor (n-3)/3\rfloor\}=\min\{\level(\la),\lfloor n/3\rfloor\}-1. 
\end{equation}

{\em Case 1:} $\level(\la)>\lfloor n/3\rfloor$.

In particular, the right hand side of (\ref{E}) equals $\lfloor n/3\rfloor-1$. Suppose (\ref{E}) fails for all $\mu\in P$. Since $\lfloor (n-3)/3\rfloor=\lfloor n/3\rfloor-1$, we have $\level(\mu)\leq \lfloor n/3\rfloor-2\leq \level(\la)-3$. But $\level((\tilde e_i^2\la)_C)\geq \level(\la)-2$ and $\level(\tilde e_{i+1}^2(\la_A)\geq \level(\la)-2$, so we get a contradiction with (\ref{EThx1}) and (\ref{EThx2}).

{\em Case 2:} $\level(\la)\leq\lfloor n/3\rfloor$. 
The right hand side of (\ref{E}) equals $\level(\la)$, and we have to find 
$\mu\in P$ such that 
\begin{equation}\label{E291022}
\min\{\level(\mu),\lfloor (n-3)/3\rfloor\}=\level(\la)-1.
\end{equation} 

If $\la_1-\la_2\leq 2$ then, since $\la_1=n-\level(\la)$, we have 
$$n\geq \la_1+\la_2\geq 2\la_1-2=2(n-\level(\la))-2\geq
2(n-\lfloor n/3\rfloor)-2
$$
which leads to the cases $\la=(3,1)$ or $(4,2)$. For $\la=(3,1)$ the equality (\ref{E291022}) is clear. For $\la=(4,2)$, using \cite[Tables]{JamesBook}, it is easy to see that $D^{(2,1)}$ is a composition factor of $e_0^2e_1D^{(4,2)}$, so in view of Lemma~\ref{LBr2}, we have $(2,1)\in P$, and so (\ref{E291022}) hols. Thus we may assume that $\la_1-\la_2\leq 3$. 

Now, if $j\neq i$ then $B$ and $C$ are $j$-normal nodes of $\la_A$. Let $D$ be the $j$-good node of $\la_A$. Then by Lemma~\ref{Lemma39}, we have that $D^{\la_{A,C,D}}$ is a composition factor of $e_j^2e_iD^\la$, so by Lemma~\ref{LBr2}, we have that $\la_{A,C,D}\in P$. But $\la_{A,C,D}=\level(\la)-1$ and (\ref{E291022}) follows. 

Finally, let $j=i$. Let $D$ be the $i$-good node of $\la$. Then by Lemma~\ref{Lemma39}, we have that $D^{\la_{A,C,D}}$ is a composition factor of $e_{i+1}e_i^2D^\la$, so by Lemma~\ref{LBr2}, we have that $\la_{A,C,D}\in P$. But $\level(\la_{A,C,D})=\level(\la)-1$ and (\ref{E291022}) follows. 
\end{proof}

\subsection{A new dimension bound}\label{SSBound}
Throughout the subsection we assume that $p>0$. Let $\la\in\Parp(n)$ and set $r:=\level(\la)$. 
%We denote by $\la^\Mull$ the Mullineux image of $\la$ so that $D^{\la^\Mull}\cong D^\la\otimes\sgn$, and $\la^\Mull=\la$ of $p=2$. The partition $\la^\Mull$ is in principle known, see \cite{FK}. Denote $r:=\level(\la)$ and $s:=\level(\la^\Mull)$. 
Then by \cite[Theorem 5.1]{GLT1}, for $n\geq 5$, we have 
$$\dim D^\la \geq 2^{r/2}.$$ 
This bound and its slight improvements \cite[Theorems B, C]{KMT} work in all characteristics, but they are far from optimal. 

Under the {\em extra assumption}\, that $p > 2$ and $n \geq p(r-2)$, or $p=2$ and $n \geq 2(r-1)$, 
\cite[Theorem A]{KMT} gives a much better lower bound 
$$\dim D^\la \geq 
\left\{
\begin{array}{ll}
\frac{n(n-p) \ldots (n-(r-1)p)}{r!} &\hbox{if $p>2$,}
\vspace{2mm}
\\
\frac{(n-p)(n-2p) \ldots (n-rp))}{r!} &\hbox{if $p=2$.}
\end{array}
\right.
$$
which is asymptotically best possible. However, the above extra assumption on $n$ limits its use, especially when 
$p > 2$.

Combining Theorems \ref{TDimRank}, \ref{TRL}, and \ref{TRL2}, we obtain a new dimension bound, whose strength is in between 
the above two bounds, but no restrictions on $n$ are needed:
%can be applied to a much wider range of irreducible modules than \cite[Theorem A]{KMT} when $p >2$:

\begin{Theorem}\label{TBound}
Let $\la\in\Parp(n)$ and $r=\level(\la)$. Then
$$
\dim D^\la \geq
\left\{
\begin{array}{ll}
\displaystyle\binom{\lfloor n/2 \rfloor}{r} &\hbox{if $p>2$,}
\vspace{2mm}
\\
\displaystyle2^r\binom{\lfloor n/3 \rfloor}{r} &\hbox{if $p=2$.}
\end{array}
\right.
$$
\iffalse{
\begin{enumerate}
\item[{\rm (i)}] If $p \neq 2$ %and $\level(\la) \leq n/2$ 
then
$$\dim D^\la \geq \binom{\lfloor n/2 \rfloor}{\level(\la)}.$$
\item[{\rm (ii)}] If $p = 2$ %and $\level(\la) \leq n/3$ 
then
$$\dim D^\la \geq 2^r\binom{\lfloor n/3 \rfloor}{\level(\la)}.$$
\end{enumerate}
}\fi 
\end{Theorem}

 \vspace{2mm}

\section{Modular analogue of the Murnaghan-Littlewood Theorem}
Theorem~\ref{TTensRank} establishes the `additivity' property of ranks for tensor products (assuming that the ranks are not too large compared to $n$). In view of Proposition~\ref{PSpechtRank} (in characteristic $0$) and Theorems~\ref{TRL},\,\ref{TRL2} (in positive characteristic), this can be restated as the 
`additivity' property of levels for tensor products of irreducible modules. In this section we substantially strengthen this additivity property proving a modular analogue of Muraghan-Littlewood theorem for tensor products of irreducible representations of symmetric groups. This must be of independent interest.

For a partition $\la=(\la_1,\la_2,\dots,\la_h)\in\Par(n)$ we define the partition 
\begin{equation}\label{ELaBar}
 \bar\la:=(\la_2,\dots,\la_h)\in\Par(\level(\la)).
\end{equation}
In this section we prove the following 

 \begin{Theorem} \label{TMur}
 Let $\al,\be,\ga\in\Par_p(n)$ satisfy $\level(\ga)=\level(\al)+\level(\be)$. Then 
 $$
 [D^\al\otimes D^\be: D^\ga]=[\Ind_{\Si_{\level(\al)}\times \Si_{\level(\be)}}^{\Si_{\level(\ga)}} (D^{\bar\al}\boxtimes D^{\bar\be}): D^{\bar\ga}].
 $$
 \end{Theorem}
 
 In characteristic $0$ this is well known and goes back to Murnaghan, see \cite{Mu} where it is stated without proof. For the first published proof see Littlewood \cite{Li}. For other proofs see for example Dvir \cite[Theorem 3.3]{Dvir} and Brion \cite[Corollary 2 in \S3.4]{Brion}. Brion's proof uses algebraic groups, and our approach to Theorem~\ref{TMur} is inspired by it.
 
 \subsection{General linear groups}
 \label{SSGL}
Let $V=\F^M$ be a finite dimensional vector space of dimension $M$ with standard basis $v_1,\dots,v_M$. We identify the groups $\GL(V)$ and $\GL_M(\F)$ using this standard basis. 
The subgroup 
$$
\T_V:=\{\diag(t_1,\dots,t_M)\mid t_1,\dots,t_M\in\F^\times\}\leq \GL(V)
$$ 
is then our standard choice of a maximal torus in $\GL(V)$.

Let $W=\F^N$ be another finite dimensional vector space with standard basis $w_1,\dots,w_N$. We again identify the groups $\GL(W)$ and $\GL_N(\F)$ using the standard basis. Moreover, we identify 
$\GL(V\oplus W)$ with $\GL_{M+N}(\F)$ using the basis 
$$v_1,\dots,v_M,w_1,\dots,w_N$$ of $V\oplus W$, and $\GL(V\otimes W)$ with $\GL_{MN}(\F)$ using the basis 
\begin{equation}\label{ETensBasis}
\{v_i\otimes w_j\mid 1\leq i\leq M,\ 1\leq j\leq N\}
\end{equation}
 of $V\otimes W$ {\em taken in lexicographic order}. 
Then the standard homomorphisms  
\begin{align*}
\iota^{\oplus}&:\GL(V)\times \GL(W)\,\into\,\GL(V\oplus W),
\\
\iota^{\otimes}&:\GL(V)\times \GL(W)\,\into\,\GL(V\otimes W).
\end{align*}
yield the restriction functors
\begin{align*}
\dar^{\GL(V\oplus W)}_{\GL(V)\times \GL(W)}&:\mod{\GL(V\oplus W)}\to \mod{\GL(V)\times \GL(W)},
\\
\dar^{\GL(V\otimes W)}_{\GL(V)\times \GL(W)}&:\mod{\GL(V\otimes W)}\to \mod{\GL(V)\times \GL(W)}.
\end{align*}
The homomorphisms $\iota^{\oplus}$ and $\iota^{\otimes}$ 
restrict to the homomorphisms (denoted by the same symbols) of the maximal tori:
\begin{align*}
\iota^{\oplus}:\T_V\times \T_W\,\iso\,\T_{V\oplus W},\ &(\diag(t_1,\dots,t_M), \diag(s_1,\dots,s_N))
\\ &\mapsto\, \diag(t_1,\dots,t_M,s_1,\dots,s_N),
\\
\iota^{\otimes}:\T_V\times \T_W\,\into\,\T_{V\otimes W},\ &(\diag(t_1,\dots,t_M), \diag(s_1,\dots,s_N))
\\ &\mapsto\, \diag(t_is_j\mid 1\leq i\leq M,\ 1\leq j\leq N),
\end{align*}
where $\iota^{\otimes}$ depends on the choice of order on the basis (\ref{ETensBasis}). 

%For a $\GL(V)$-module $X$ and a $\GL(W)$-module $Y$ we denote by $X\boxtimes Y$ the $\GL(V)\times\GL(W)$-module which is the outer tensor product of $X$ and $Y$. 

We will consider polynomial weights of maximal tori and identify them with compositions as follows. Let $\La(M,m)$ be the set of 
all compositions $\mu=(\mu_1,\dots,\mu_M)\in\Z_{\geq 0}^M$ such that $\mu_1+\dots+\mu_M=m$, see \cite[\S3]{Green}. Then we identify $\mu$ with the degree $m$ polynomial weight 
$$
\mu:\T_V\to\F^\times,\ \diag(t_1,\dots,t_M)\mapsto t_1^{\mu_1}\cdots t_M^{\mu_M}
$$
of $\T_V$.  
Similarly we identify degree $n$ polynomial weights of $\T_W$ with the set $\La(N,n)$. Finally, we identify the degree $k$ polynomial weights of $\T_{V\otimes W}$ with the set $\La(M\times N,k)$ of tuples $\la=(\la_{i,j})_{\substack{1\leq i\leq M\\ 1\leq j\leq N}}$ of non-negative integers summing to $n$ as follows: given $\la\in \La(M\times N,k)$ we identify it with the weight
$$
\la:\T_{V\otimes W}\to\F^\times,\ \diag(t_{i,j}\mid 1\leq i\leq M,\ 1\leq j\leq N)\mapsto \prod_{i=1}^M\prod_{j=1}^Nt_{i,j}^{\la_{i,j}}.
$$

Let $\mu\in\La(M,m)$ and $\nu\in\La(N,n)$. We have the weight 
$$\mu\times\nu:=(\mu_1,\dots,\mu_M,\nu_1,\dots,\nu_N)\in\La(M+N,m+n)$$
of $T_{V\oplus W}=T_V\times T_W$ 
On the other hand, inflating $\la\in\La(M\times N,k)$ along the embedding $\iota^\otimes: \T_V\times \T_W\into \T_{V\otimes W}$ gives the polynomial weight $\infl^\otimes \la$. Note that 
$$
\infl^\otimes \la:=\mu\times \nu\in\La(M+N,2k),
$$
where
\begin{equation}\label{EInflF}
\mu_i=\sum_{j=1}^N\la_{ij}\ \text{for all $1\leq i\leq M$ \qquad and}\qquad\   
\nu_j=\sum_{i=1}^M\la_{ij}\ \text{for all $1\leq j\leq N$}.
\end{equation}
%If these conditions hold, we will write $\la\in\Wt_{\mu,\nu}(V\otimes W)$. 

\subsection{Polynomial modules over general linear groups}
For $m\in\Z_{\geq 0}$ and $V=\F^M$ as in \S\ref{SSGL}, 
we consider the category $\Pol(V,m)$ of degree $m$ polynomial representations of $\GL(V)$, see \cite[\S2.2]{Green}. The set of weights of modules in $\Pol(V,m)$ lies in $\La(M,m)$, see \cite[\S3]{Green}. For $X\in \Pol(V,m)$ and $\la\in\La(M,m)$, we have the $\la$-weight space 
$$X_\la:=\{x\in X\mid \diag(t_1,\dots,t_M)x=t_1^{\la_1}\cdots\la_M^{\la_M}x\ \text{for all}\  t_1,\dots,t_n\in\F^\times\}.
$$
Then $X=\bigoplus_{\la\in\La(M,m)}X_\la$. The {\em formal character} of $X$ is defined as
\begin{equation}\label{ECh}
\ch X:=\sum_{\la\in \La(M,m)}(\dim X_\la)e^\la.
\end{equation}

We will often suppose that $m\leq M$. In this case the subset of dominant weights $\La^+(M,m)\subseteq \La(M,m)$ as in \cite[\S3.1]{Green} can be identified with the set $\Par(m)$ of partitions of $m$ via the bijection sending the partition $\al=(\al_1\geq\dots\geq\al_h>0)\in \Par(m)$ to $(\al_1,\dots,\al_h,0,\dots,0)\in\La^+(M,m)$. 

In general, for any $\al\in \La^+(M,m)$, 
%Thus for $m\leq M$ and $\al\in\Par(m)$, 
we have the following canonical modules in $\Pol(V,m)$:
\begin{enumerate}
\item[$\bullet$] $L_\al(V)$, the irreducible module with highest weight $\al$;
\item[$\bullet$] $\De_\al(V)$, the standard module with highest weight $\al$;
\item[$\bullet$] $\nabla_\al(V)$, the costandard module with highest weight $\al$;
\item[$\bullet$] $T_\al(V)$, the indecomposable tilting module with highest weight $\al$.
\end{enumerate}
For irreducible, standard and costandard modules, see \cite{Green}. For tilting modules see \cite{Donkin}. 

We can express the multiplicity in the right hand side of the equality in Theorem~\ref{TMur} in terms of tilting modules using the following lemma:

\begin{Lemma} \label{LBKMLR} {\rm \cite[Theorem D]{BKLR}}
Let $\al\in\Par_p(m)$, $\be\in\Par_p(n)$, $\ga\in\Par_p(m+n)$, $V=\F^M$ and $W=\F^N$ such that $m\leq M$ and $n\leq N$. Then the restriction $T_\ga(V\oplus W)\dar^{\GL(V\oplus W)}_{\GL(V)\times\GL(W)}$ is tilting, and 
$$
\big(T_\ga(V\oplus W)\dar^{\GL(V\oplus W)}_{\GL(V)\times\GL(W)}:T_\al(V)\boxtimes T_\be(W)\big)=\big[\Ind_{\Si_{m}\times\Si_n}^{\Si_{m+n}}(D^\al\boxtimes D^\be):D^\ga \big].
$$
\end{Lemma}

The multiplicity in the left hand side of Lemma~\ref{LBKMLR} (even without assuming that $\al,\be,\ga$ are $p$-regular) can be related to another ``modular Littlewood-Richardson coefficient":

\begin{Lemma} \label{L4'} {\rm \cite[Theorem B(i)]{BKLR}}
Let $\al\in\Par(m)$, $\be\in\Par(n)$, $\ga\in\Par(m+n)$, $V=\F^M$,  $W=\F^N$ and $U=\F^K$ such that $m\leq M$, $n\leq N$ and $m+n\leq K$. Then the restriction $T_\ga(V\oplus W)\dar^{\GL(V\oplus W)}_{\GL(V)\times\GL(W)}$ is tilting, and 
$$
\big(T_\ga(V\oplus W)\dar^{\GL(V\oplus W)}_{\GL(V)\times\GL(W)}:T_\al(V)\boxtimes T_\be(W)\big)=\big[L_{\al'}(U)\otimes L_{\be'}(U):L_{\ga'}(U)].
$$
\end{Lemma}

Weight multiplicities of tilting modules can also be expressed in terms of certain ``modular Littlewood-Richardson coefficients". For a composition $\mu=(\mu_1,\dots,\mu_M)\in\La(M,m)$ and a vector space $U$ we denoted 
$$
{\textstyle\bigwedge}^\mu U:={\textstyle\bigwedge}^{\mu_1}U\otimes\dots\otimes {\textstyle\bigwedge}^{\mu_N}U,
$$
considered as a module in $\Pol(U,m)$. Then:

\begin{Lemma} \label{L4} {\rm \cite[Theorem C(iii)]{BKLR}}
Let $\al\in\Par(m)$, $\mu\in\La(M,m)$, $V=\F^M$ and $U=\F^K$ such that $m\leq M,K$. Then 
$$
\dim T_\al(V)_\mu=[{\textstyle\bigwedge}^\mu U:L_{\al'}(U)].
$$
\end{Lemma}

For weight multiplicities of tilting modules we also have: 

\begin{Lemma} \label{L3'} {\rm \cite[(1.5)(ii)]{Donkin}}
Let $\al\in\La^+(M,m)$ and $\mu\in\La(M,m)$ be such that $\al_1=\mu_1$. Let $V=\F^M$, $\tilde V=\F^{M-1}$ Considering $\tilde\al:=(\al_2,\dots,\al_M)\in\La^+(M-1,m-\al_1)$ and 
$\tilde\mu:=(\mu_2,\dots,\mu_M)\in\La(M-1,m-\al_1)$ as weights for $\GL(\tilde V)$, we have 
$$
\dim T_\al(V)_\mu=\dim T_{\tilde\al}(\tilde V)_{\tilde\mu}.
$$
\end{Lemma}

\begin{Corollary} \label{C3} 
Let $m\leq N$, $\al\in\La^+(M,m)$ and $\mu\in\La(M,m)$ be such that $\al_1=\mu_1$. Let $V=\F^M$. Considering $\bar\al:=(\al_2,\dots,\al_M,0)\in\La^+(M,m-\al_1)$ and 
$\bar\mu:=(\mu_2,\dots,\mu_M, 0)\in\La(M,m-\al_1)$ as weights for $\GL(V)$, we have 
$$
\dim T_\al(V)_\mu=\dim T_{\bar\al}(V)_{\bar\mu}.
$$
\end{Corollary}
\begin{proof}
By Lemma~\ref{L3'}, we have 
$
\dim T_\al(V)_\mu=\dim T_{\tilde\al}(V)_{\tilde\mu}.
$
where $\tilde\al=(\al_2,\dots,\al_M)$ and $\tilde\mu=(\mu_2,\dots,\mu_M)$ are considered as weights of $GL(\tilde V)$. It remains to note that 
$\dim T_{\tilde\al}(\tilde V)_{\tilde\mu}=\dim T_{\bar\al}(V)_{\bar\mu}
$
for example using Lemma~\ref{L4}.
\end{proof}

\subsection{Schur-Weyl duality}
We now consider $V^{\otimes m}\in\Pol(V,m)$. Let $m\leq M=\dim V$. Then the endomorphism algebra $\End_{\F\GL(V)}(V^{\otimes m})$ is naturally identified with the group algebra $\F\Si_m$ (via the action of $\Si_m$ on tensors by permuting components), see e.g. \cite{Green}. 

Since the natural module $V$ is a tilting module, $V^{\otimes m}$ is also tilting. Moreover, by \cite[Proposition 4.2]{Erdmann}, we have 
\begin{equation}\label{EErdmann}
V^{\otimes m}\cong\bigoplus_{\la\in\Par_p(m)}T_\la(V)^{\bigoplus \dim D^\la},
\end{equation}
We will need to slightly improve on this result in Lemma~\ref{LErdmannImprove} below. 
To explain this, we consider the following general set up. 

Let $A$ be an $\F$-algebra and $T$ be a finite dimensional $A$-module. Let 
\begin{equation}\label{ETi}
T=\bigoplus_{i\in I}T_i
\end{equation} 
be a decomposition into a direct sum of indecomposable $A$-modules. Define an equivalence relation $\sim$ on $I$ by $i\sim j$ if and only if $T_i\cong T_j$. 

We have the centralizer algebra $$C:=\End_A(T)^\op$$ with idempotents $\{e_i\mid i\in I\}$, where, for each $i\in I$, we define  $e_i$ to be the projection onto $T_i$ along the decomposition (\ref{ETi}). 
Note that, for each $i\in I$, 
$$
P_i:=\Hom_A(M,T_i)
$$
is a left $C$-module via $c\cdot\theta=\theta\circ c$ for all $\theta\in P_i$ and $c\in C$. Clearly we have an isomorphism of $C$-modules $P_i\cong Ce_i$. 

By Fitting's Lemma, see e.g. \cite[Theorem 1.4]{Landrock}, we have that $1_C=\sum_{i\in I}e_i$ is an orthogonal decomposition of $1_C$ into primitive idempotents, and $i\sim j$ if and only if $P_i\cong P_j$ as $C$-modules. Let 
$$X:=I/\sim$$ and for $\la\in X$, let 
$$T_A(\la):=T_i\quad\text{and}\quad P_C(\la):=P_i$$ 
for some $i\in\la$. Also, let $L_C(\la)$ be the (irreducible) head of $P_C(\la)$. Then $\dim L_C(\la)=|\la|$, and 
$$\{L_C(\la)\mid \la\in X\}$$ 
is a complete set of non-isomorphic irreducible $C$-modules. So  we have an $A$-module decomposition
$$T\cong\bigoplus_{\la\in X}T_A(\la)^{\oplus \dim L_C(\la)}.$$ 

Returning to our example where $A=\F\GL(V)$, $T=V^{\otimes m}$ with $m\leq M=\dim V$, and $C=\F\Si_m$, in view of (\ref{EErdmann}), we can identify $X$ with $\Par_p(m)$ so that $T_A(\la)=T_\la(V)$ for all $\la\in\Par_p(m)$. What is not immediately clear is that then $L_C(\la)$ is identified with $D^\la$ (rather than, for example, with $D^\la\otimes\sgn$). This is established in the following lemma:

\begin{Lemma} \label{LErdmannImprove}%{\rm \cite{}}%{\bf ()}
In the notation above, we have $L_C(\la)\cong D^\la$. 
\end{Lemma}
\begin{proof}
We have identified $\Par_p(m)$ with $X$ so that  $e_\la\in\F\Si_m=\End_{\F\GL(V)}(V^{\otimes m})$ is a projection onto  the indecomposable summand isomorphic to $T_\la(V)$. To show that $L_C(\la)\cong D^\la$, it suffices to prove that the $\F\Si_m$-module 
$$
P_C(\la)=\Hom_{\F\GL(V)}(V^{\otimes d},T_\la(V))
$$
contains $D^\la$ in its head. 

Since $T_\la(V)$ has a $\nabla$-filtration and $\la$ is its highest weight, we can find a submodule $K\subset T_n(\la)$ with a $\nabla$-filtration such that we have an exact sequence
$$
0\to K\to T_\la(V)\to \nabla_\la(V)\to 0.
$$
Since $V^{\otimes m}$ has a $\De$-filtration, we have $\Ext^1_{\Pol(V,m)}(V^{\otimes m}, K)=0$. So, 
applying the functor $\Hom_{\F\GL(V)}(V^{\otimes m},-)$ to this exact sequence, we obtain an exact sequence of $\F\Si_m$-modules 
$$
\Hom_{\F\GL(V)}(V^{\otimes m},T_\la(V))\to \Hom_{\F\GL(V)}(V^{\otimes m},\nabla_\la(V))\to 0.
$$
But it is well known, see e.g. \cite[\S6.3]{Green}, that the $\F\Si_m$-module $\Hom_{\F\GL(V)}(V^{\otimes m},\nabla_\la(V))$ is isomorphic to the Specht module $S^\la$ so it has $D^\la$ in its head. 
\end{proof}

Returning to the general situation, let now $B$ be a subalgebra of $A$. Then $C$ is naturally a subalgebra of the endomorphism algebra 
$$
D:=\End_B(T)^\op.
$$
Applying the theory developed above to the pair $(B,D)$ instead of the pair $(A,C)$, we get the set of equivalence classes $Y$, the indecomposable modules $T_B(\mu)$ and the irreducible modules $L_D(\mu)$ with 
$$T\cong\bigoplus_{\mu\in Y}T_B(\mu)^{\oplus \dim L_D(\mu)}.$$ 

The following useful general observation was suggested by O. Mathieu \cite{Mathieu}:

\begin{Lemma} \label{LMathieu} {\rm \cite[Lemma 3.1]{BKLR}} 
For $\la\in X$, we have 
$$
T_A(\la)\dar_B\cong\bigoplus_{\mu\in Y}T_B(\mu)^{\oplus [L_D(\mu)\dar_C:L_C(\la)]}.
$$
In other words, $(T_A(\la)\dar_B:T_B(\mu))=[L_D(\mu)\dar_C:L_C(\la)]$ for all $\mu\in Y$. 
\end{Lemma}

This lemma will be applied in the following situation. Let $V=\F^M$ and $W=\F^N$. Take:
\begin{enumerate}
\item[$\bullet$] $A$ is the group algebra $\F\GL(V\otimes W)$,
\item[$\bullet$] $T=(V\otimes W)^{\otimes n}$,
\item[$\bullet$] $B$ is the image of the group algebra $\F(\GL(V)\times \GL(W))$ under the homomorphism $\iota^\otimes$ from $\F(\GL(V)\times \GL(W))$ to $\F\GL(V\otimes W)$. 
 \end{enumerate} 
We will consider $B$-modules as $\F(\GL(V)\times \GL(W))$-modules via inflation. 
We assume that $n\leq M,N$. Then certainly $n\leq MN$, so $C=\F\Si_n$. On the other hand, we have
\begin{align*}
D&=\End_{\GL(V)\times GL(W)}((V\otimes W)^{\otimes n})
\\
&\cong \End_{\GL(V)\times GL(W)}(V^{\otimes n}\boxtimes W^{\otimes n})
\\
&\cong \End_{\GL(V)}(V^{\otimes n})\otimes 
\End_{\GL(W)}(W^{\otimes n})
\\
&\cong \F\Si_n\otimes 
\F\Si_n,
\end{align*}
and the natural embedding of $C=\F\Si_n$ into $D=\F\Si_n\otimes\F\Si_n$ comes from the diagonal embedding $\Si_n\to\Si_n\times\Si_n$. So given $S,T\in\mod{\F\Si_n}$, restricting $S\boxtimes T$, from $D$ to $C$ yields the inner product $S\otimes T$. Now Lemmas~\ref{LMathieu} and \ref{LErdmannImprove} imply:

\begin{Theorem} \label{TBrion}
Let $\al,\be,\ga\in\Par_p(n)$, $V=\F^M$ and $W=\F^N$ with $M,N\geq n$. Then the restriction $T_\ga(V\otimes W)\dar^{\GL(V\otimes W)}_{\GL(V)\times\GL(W)}$ is tilting, and
$$
\big(T_\ga(V\otimes W)\dar^{\GL(V\otimes W)}_{\GL(V)\times \GL(W)}:T_\al(V)\boxtimes T_\be(W)\big)=\big[D^\al\otimes D^\be:D^\ga\big]. 
$$
\end{Theorem}

Note that the partition $\ga$ in the theorem is considered as a dominant weight for the torus $\T_{V\otimes W}$, see \S\ref{SSGL}. 

\subsection{Lemma on weights}
Throughout the subsection, $V=\F^M$ and $W=\F^N$, $n\leq M,N$, and 
$\ga$ is a partition of $n$. Theorem~\ref{TBrion} explains out interest in the restriction in 
$$T_\ga(V\otimes W)\dar^{\GL(V\otimes W)}_{\GL(V)\times \GL(W)}.
$$  

In this subsection we analyze weights of this restriction. 
%Recall that $\mu=(\mu_1,\dots,\mu_M)\in\La(M,k)$ is considered as a degree $k$ polynomial weight of $\T_V$ and $\nu=(\nu_1,\dots,\nu_N)\in\La(N,l)$ is considered as a degree $l$ polynomial weight of $\T_W$, so $$\mu\times\nu:=(\mu_1,\dots,\mu_M,\nu_1,\dots,\nu_N)\in\La(M+N,k+l)$$ can be considered as a weight of $T_{V\oplus W}=T_V\times T_W$ in the obvious way. 

%On the other hand, degree $n$ polynomial weights of the maximal torus $T_{V\otimes W}$ can be identified in the obvious way with tuples $\la=(\la_{i,j})_{\substack{1\leq i\leq M\\ 1\leq j\leq N}}$ of non-negative integers summing to $n$. We denote the set of such weights by $\La(V\otimes W,n)$. 
%Inflating $\la\in\La(V\otimes W,n)$ along the embedding $\iota^\otimes: \T_V\times \T_W\into \T_{V\otimes W}$ gives the polynomial weight $\infl^\otimes \la$. Note that 
%$$\infl^\otimes \la:=\mu\times \nu\in\La(M+N,2n),$$
%where
%$$\mu_i=\sum_{j=1}^N\la_{ij}\ \text{for all $1\leq i\leq M$ \qquad and}\qquad\   \nu_j=\sum_{i=1}^M\la_{ij}\ \text{for all $1\leq j\leq N$}.$$

\begin{Lemma} \label{LWeights}
Let $M,N\geq n$, $\ga\in\Par(n)$,  
$\mu\in\La(M,n)$ and $\nu\in\La(N,n)$ with 
$a:=\mu_2+\dots+\mu_M$ and $b:=\nu_2+\dots+\nu_N$ satisfying 
$a+b\leq \level(\ga)$. If $\mu\times\nu=\infl^\otimes\la$ for some weight $\la$ of $T_\ga(V\otimes W)$ then:
\begin{enumerate}
\item[{\rm (1)}] $a+b=\level(\ga)$. 
\item[{\rm (2)}] $\la_{1,1}=n-a-b$,
\item[{\rm (3)}] $\la_{1,j}=\nu_j$ for $j=2,\dots,N$,
\item[{\rm (4)}] $\la_{i,1}=\mu_i$ for $i=2,\dots,M$,
\item[{\rm (5)}] $\la_{i,j}=0$ for $i,j\geq 2$.
\end{enumerate}
\end{Lemma}
\begin{proof}
The highest weight of $T_\ga(V\otimes W)$ is $\ga$ 
(considered as a weight of $\T_{V\otimes W}$). 
Since $\ga_1=n-\level(\ga)$ and the dominant conjugate of the weight $\la$ is less than or equal to $\ga$ in the dominance order, it follows that $\la_{i,j}\leq n-\level(\ga)$ for all $i,j$. 
On the other hand, recalling (\ref{EInflF}), we have 
\begin{align*}
n-a&=\mu_1
\\
&=\la_{1,1}+\la_{1,2}+\dots+\la_{1,N}
\\
&\leq n-\level(\ga)+\la_{1,2}+\dots+\la_{1,N}
\\
&\leq n-a-b+\la_{1,2}+\dots+\la_{1,N}
\\
&\leq n-a-b+\sum_{i=1}^M\la_{i,2}+\dots+\sum_{i=1}^M\la_{i,N}
\\
&= n-a-b+\nu_{2}+\dots+\nu_{N}
\\
&= n-a-b+b
\\
&= n-a.
\end{align*}
This proves that $a+b=\level(\ga)$, $\la_{1,1}= n-\level(\ga)= n-a-b$, and $\la_{1,j}=\sum_{i=1}^M\la_{i,j}=\nu_j$ for $j=2,\dots,N$; in particular, we get (1), (2), (3) and (5). The claim (4) follows from (5) and (\ref{EInflF}).
\end{proof}

\subsection{Proof of Theorem~\ref{TMur}} 
Throughout the subsection, $V=\F^M$ and $W=\F^N$ with $M,N\geq n$. We fix $p$-regular partitions $\al,\be,\ga$ of $n$. 
For $a\in\Z_{\geq 0}$ we will use the notation
\begin{align*}
\La(M,n;a)&:=\{\mu\in \La(M,n)\mid \mu_2+\dots+\mu_M=a\},
\\
\La(M,n;\leq a)&:=\{\mu\in \La(M,n)\mid \mu_2+\dots+\mu_M\leq a\},
\\
\La^+(M,n;a)&:=\La(M,n;a)\cap\La^+(M,n),
\\
\La^+(M,n;\leq a)&:=\La(M,n;\leq a)\cap\La^+(M,n).
\end{align*}
The sets $\La(N,n;a),\La(N,n;\leq a)$ etc. are defined similarly.

By Theorem~\ref{TBrion}, the left hand side of the equality in Theorem~\ref{TMur} equals 
$$
 \big(T_\ga(V\otimes W)\dar^{\GL(V\otimes W)}_{\GL(V)\times \GL(W)}:T_\al(V)\boxtimes T_\be(W)\big).
$$
On the other hand, by Lemma~\ref{LBKMLR}, the right hand side of the equality in Theorem~\ref{TMur} equals 
$$
\big(T_{\bar\ga}(V\oplus W)\dar^{\GL(V\oplus W)}_{\GL(V)\times\GL(W)}:T_{\bar\al}(V)\boxtimes T_{\bar\be}(W)\big).
$$
So Theorem~\ref{TMur} follows from the following result on tilting modules:

\begin{Theorem} \label{TTilt}
 Let $\al,\be,\ga\in\Par_p(n)$ satisfy $\level(\ga)=\level(\al)+\level(\be)$. Then 
 $$
 \big(T_\ga(V\otimes W)\dar^{\GL(V\otimes W)}_{\GL(V)\times \GL(W)}:T_\al(V)\boxtimes T_\be(W)\big)
 =\big(T_{\bar\ga}(V\oplus W)\dar^{\GL(V\oplus W)}_{\GL(V)\times\GL(W)}:T_{\bar\al}(V)\boxtimes T_{\bar\be}(W)\big).
$$
\end{Theorem}
\begin{proof}
Set $a:=\level(\al)$ and $b:=\level(\be)$. For any polynomial dominant weights $\theta$ of $\T_V$ and $\eta$ of $\T_W$, we denote
$$
m^\ga_{\theta,\eta}:=\big(T_\ga(V\otimes W)\dar^{\GL(V\otimes W)}_{\GL(V)\times \GL(W)}:T_\theta(V)\boxtimes T_\eta(W)\big).
$$
We can write
$$
T_\ga(V\otimes W)\dar^{\GL(V\otimes W)}_{\GL(V)\times \GL(W)}=
\bigoplus_{\substack{\theta\in\La^+(M,n;\leq a)\\ \eta\in\La^+(N,n;\leq b)}}
\big(T_\theta(V)\boxtimes T_\eta(W)\big)^{\oplus m^\ga_{\theta,\eta}}\,\oplus\,X,
$$
where $X$ stands for the sum of the remaining summands  $T_\theta(V)\boxtimes T_\eta(W)$ of $T_\ga(V\otimes W)\dar^{\GL(V\otimes W)}_{\GL(V)\times \GL(W)}$. For $\mu\in\La(M,n;\leq a)$ and $\nu\in\La(N,n;\leq b)$ we then have $X_{\mu\times\nu}=0$. So, we have for the dimensions of $(\mu\times\nu)$-weight spaces:
\begin{equation}\label{E*}
\dim \Big(T_\ga(V\otimes W)\dar^{\GL(V\otimes W)}_{\GL(V)\times \GL(W)}\Big)_{\mu\times\nu}=
\sum_{\substack{\theta\in\La^+(M,n;\leq a)\\ \eta\in\La^+(N,n;\leq b)}}
m^\ga_{\theta,\eta}
\dim T_\theta(V)_\mu\,\dim T_\eta(W)_\nu.
\end{equation}

Recalling (\ref{ECh}), we now consider the truncated characters  
$$
\overline{\ch} X:=\sum_{\mu\in\La(M,n;\leq a)}(\dim X_\mu)e^\mu\quad\text{and}\quad 
\overline{\ch} Y:=\sum_{\nu\in\La(N,n;\leq b)}(\dim X_\nu)e^\nu
$$ 
for $X\in\Pol(V,n)$ and $Y\in \Pol(W,n)$. Since the highest weight of $T_\theta(V)$ is $\theta$, the truncated characters
$$
\{\overline{\ch} \,T_\theta(V)\mid\theta\in\La^+(M,n;\leq a)\}
$$
are linearly independent. Similarly the truncated characters
$$
\{\overline{\ch} \,T_\eta(W)\mid\theta\in\La^+(N,n;\leq b)\}
$$
are linearly independent. It follows that the equations (\ref{E*}) for all $\mu\in\La(M,n;\leq a)$ and $\nu\in\La(N,n;\leq b)$ determine the coefficients $m^\la_{\theta,\eta}$. 

Now, $\mu\times\nu$ is a weight of $T_\ga(V\otimes W)\dar^{\GL(V\otimes W)}_{\GL(V)\times \GL(W)}$ if and only if $\mu\times\nu=\infl^\otimes \la$ for some weight $\la$ of $T_\ga(V\otimes W)$. So for $\mu\in\La(M,n;\leq a)$ and $\nu\in\La(N,n;\leq b)$, by Lemma~\ref{LWeights}, the weight space in the left hand side of (\ref{E*}) is non-trivial only if
\begin{enumerate}
\item[{\sf (a)}] $\mu_2+\dots+\mu_M+\nu_2+\dots+\nu_N=\level(\ga)$.
\item[{\sf (b)}] $\la$ is uniquely determined by $\mu$ and $\nu$ via the equalities (2)-(5) of Lemma~\ref{LWeights}. We denote such $\la$ by $\la(\mu,\nu)$. Then $\la(\la,\mu)_{1,1}=\ga_1$, $\la(\mu,\nu)_{1,j}=\nu_j$ for $j=2,\dots,N$, $\la(\mu,\nu)_{i,1}=\mu_i$  for $i=2,\dots,M$, and $\la(\mu,\nu)_{i,j}=0$ otherwise.
\end{enumerate}
But $\mu_2+\dots+\mu_M\leq a$ since $\mu\in\La(M,n;\leq a)$ and similarly $\nu_2+\dots+\nu_N\leq b$. On the other hand, $a+b= \level(\ga)$ by assumption. So we conclude from (a) that $\mu\in\La(M,n; a)$ and $\nu\in\La(N,n;b)$. 
 
Since $T_\theta(V)_\theta\neq 0\neq T_\eta(W)_\eta$ we deduce that $m^\la_{\theta,\eta}=0$ for $\theta\in\La^+(M,n;\leq a)$ and $\eta\in\La^+(N,n;\leq b)$ unless 
$\theta\in\La^+(M,n;a)$ and $\eta\in\La^+(N,n;b)$. Moreover, in view of {\sf (b)}, we now get from (\ref{E*}).

\begin{equation}\label{E**}
\dim T_\ga(V\otimes W)_{\la(\mu,\nu)}=
\sum_{\substack{\theta\in\La^+(M,n;a)\\ \eta\in\La^+(N,n;b)}}
m^\ga_{\theta,\eta}
\dim T_\theta(V)_\mu\,\dim T_\eta(W)_\nu.
\end{equation}
Again, by linear independence of truncated characters, the equations (\ref{E**}) for all $\mu\in\La(M,n; a)$ and $\nu\in\La(N,n;b)$ determine the coefficients $m^\la_{\theta,\eta}$.

For $\mu\in\La(M,n; a)$ and $\nu\in\La(N,n;b)$, denoting $\bar\mu:=(\mu_2,\dots,\mu_M)$ and $\bar\nu:=(\nu_2,\dots,\nu_N)$ we now obtain: 
\begin{align*}
&\dim T_\ga(V\otimes W)_{\la(\mu,\nu)}
\\
\stackrel{(1)}{=}
\,&\dim T_{\bar\ga}(V\otimes W)_{\bar\la(\mu,\nu)}
\\
\stackrel{(2)}{=}
\,&
\Big[\bigotimes_{i=1}^M\bigotimes_{j=1}^N{\textstyle \bigwedge^{\bar\la(\mu,\nu)_{i,j}}}\hspace{.5mm}U:L_{(\bar\ga)'}(U)\Big]
\\
\stackrel{(3)}{=}
\,&
\Big[
{\textstyle \bigwedge^{\mu_{2}}}\hspace{.5mm}U
\otimes\dots\otimes 
{\textstyle \bigwedge^{\mu_{M}}}\hspace{.5mm}U
\otimes 
{\textstyle \bigwedge^{\nu_{2}}}\hspace{.5mm}U
\otimes\dots\otimes
{\textstyle \bigwedge^{\nu_{N}}}\hspace{.5mm}U
:L_{(\bar \ga)'}(U)\Big]
\\
\stackrel{(4)}{=}
\,&
\Big[
{\textstyle \bigwedge^{\bar\mu}}\hspace{.5mm}U
\otimes 
{\textstyle \bigwedge^{\bar\nu}}\hspace{.5mm}U
:L_{(\bar \ga)'}(U)\Big].
\\
\stackrel{(5)}{=}
\,&
\sum_{\substack{\bar\theta\in\Par(a)\\ \bar\eta\in\Par(b)}}\Big[
{\textstyle \bigwedge^{\bar\mu}}\hspace{.5mm}U:L_{(\bar\theta)'}(U)\big]
\big[
{\textstyle \bigwedge^{\bar\nu}}\hspace{.5mm}U
:L_{(\bar\eta)'}(U)\big]
:[L_{(\bar\theta)'}(U)\otimes L_{(\bar\eta)'}(U):L_{(\bar \ga)'}(U)\Big]
\\
\stackrel{(6)}{=}
\,&
\sum_{\substack{\bar\theta\in\Par(a)\\ \bar\eta\in\Par(b)}}
\dim T_{\bar\theta}(V)_{\bar\mu}\,
\dim T_{\bar\eta}(V)_{\bar\nu}
\,
\Big(T_{\bar\ga}(V\oplus W)\dar^{\GL(V\oplus W)}_{\GL(V)\times\GL(W)}:T_{\bar\theta}(V)\boxtimes T_{\bar\eta}(W)\Big)
\\
\stackrel{(7)}{=}
\,&
\sum_{\substack{\theta\in\La^+(M,m;a)\\ \eta\in\La^+(N,n;b)}}
\dim T_{\theta}(V)_{\mu}\,
\dim T_{\eta}(V)_{\nu}
\,
\Big(T_{\bar\ga}(V\oplus W)\dar^{\GL(V\oplus W)}_{\GL(V)\times\GL(W)}:T_{\bar\theta}(V)\boxtimes T_{\bar\eta}(W)\Big),
\end{align*}
where:

equality (1) comes from Corollary~\ref{C3} using the equality $\la(\la,\mu)_{1,1}=\ga_1$ which we have by (b);

equality (2) comes from Lemma~\ref{L4}; 

equality (3) comes from {\sf (b)};

equality (4) comes from the notation $\bar\mu:=(\mu_2,\dots,\mu_M)$ and $\bar\nu:=(\nu_2,\dots,\nu_N)$;

equality (5) comes by considering the composition factors $L_{(\bar\theta)'}(U)$ of $\bigwedge^{\bar\mu}U$ and $L_{(\bar\eta)'}(U)$ of $\bigwedge^{\bar\nu}U$;

equality (6) comes from Lemmas~\ref{L4} and \ref{L4'};

equality (7) comes from Corollary~\ref{C3}.

Comparing to (\ref{E**}), we deduce that 
$$m^\la_{\theta,\eta}=\Big(T_{\bar\ga}(V\oplus W)\dar^{\GL(V\oplus W)}_{\GL(V)\times\GL(W)}:T_{\bar\mu}(V)\boxtimes T_{\bar\nu}(W)\Big),
$$
as required.
\end{proof}

\section{Tensor product growth of representations of symmetric and alternating groups}

Let $G$ be a finite group. Following \cite{LST}, 
the \emph{Plancherel measure} $|L|$ of an irreducible $\F G$-module  is defined to be $(\dim V^2)$.  The Plancherel measure $|V|$ of an arbitrary $\F G$-module $V$ is defined to be the sum of the Plancherel measures of its composition factors ignoring multiplicities, i.e. denoting by $\Irr(\F G)$ the set of the isomorphism classes of irreducible $\F G$-modules, we have 
$$
|V|=\sum_{L\in\Irr(\F G)\ \text{with}\ [V:L]\neq 0}|L|
=\sum_{L\in\Irr(\F G)\ \text{with}\ [V:L]\neq 0}(\dim L)^2.
$$

%$X = a_1\chi_1+\cdots+a_k\chi_k$, where the $\chi_i$ are irreducible and pairwise distinct, is defined to be $\sum_i |\chi_i|$.  Note that the multiplicities $a_i$ do not affect Plancherel measure.

\subsection{Tensor product growth of complex representations of $\Si_n$}
%Let $\lambda$ be a partition of $n$ and $\lambda'$ its dual partition.
Throughout the subsection we assume that $\F=\C$.   
Let $\chi^\lambda$ denote the character of the (irreducible) Specht module $S^\lambda$.  Recall the notation (\ref{ELaBar}). 

\begin{Lemma}
\label{lambda vs bar}
If $\lambda$ is a partition of $n$ of level $l\le n/3$, then
$$\frac 12\binom nl \le \frac{\chi^\lambda(1)}{\chi^{\bar\lambda}(1)}\le \binom nl.$$
Moreover,
\begin{equation}
\label{lower bound}
\chi^\lambda(1) \ge \Bigl(\frac{n-2l}{l}\Bigr)^{l}.
\end{equation}
\end{Lemma}

\begin{proof}
Let $P$ denote the product of all hook lengths for boxes in the first row of the Young diagram of $\lambda$. Then by the Hook Formula we have 
$$P \frac{l!}{\chi^{\bar\lambda}(1)} = \frac{n!}{\chi^\lambda(1)}.$$

Let $a_i = \lambda'_i-1$, so $a_i+1$ is the length of the $i$th column of the  Young diagram of $\la$.
We have $a_i = 0$ for $i> n-l$ and 
\begin{equation}
\label{a constraints}
\begin{split}
a_1\ge a_2 \ge \cdots\ge a_{n-l}\ge 0,\\
a_1 + a_2 + \cdots  + a_{n-l}= l,\\
a_i\in \Z.
\end{split}
\end{equation}
Let
$$P(a_1,\ldots,a_{n-l}) := (n-l+a_1-1)(n-l+a_2-2)\cdots (n-l+a_{n-l}-(n-l)).$$
We choose $(a_1,\ldots,a_{n-l})$ among the finite set of $n-l$-tuples satisfying \eqref{a constraints}
to maximize the value of the function $P$.
If $A> B-1$, then $(A+1)(B-1) < AB$, so if $a_i \ge a_{i+1}+2$, then 
we could increase the value of $P$ by decreasing $a_i$ by $1$ and increasing $a_{i+1}$ by $1$.
Therefore, $a_{i+1}\in \{a_i,a_i-1\}$.  Likewise, if $a_{i+1} = a_i-1$ and $a_{j+1} = a_j-1$ for some $0<i<j<n-l$,
then we could increase the value of $P$ by decreasing $a_i$ by $1$ and increasing $a_{j+1}$ by $1$.
Therefore, $P$ is maximized when $a_1 = \cdots = a_l = 1$ and $a_{l+1}=\cdots=a_{n-l}=0$.
This gives a value 
$$\frac{n-l+1}{n-2l+1}(n-l)! \le 2(n-l)!.$$
On the other hand, all values of $P$ are at least $(n-l)!$.  The first part of the lemma follows immediately.  

For the second part, we note that the product of all hook lengths of 
$\lambda$ is at most $(n-l+1)(n-l)\cdots (n-2l+2)\cdot (n-2l)!\cdot l!$, so
$$\chi^\lambda(1) \ge \frac{n(n-1)\cdots(n-l+2)\cdot(n-2l+1)}{l!} \ge \Bigl(\frac{n-2l}l\Bigr)^l.$$
\end{proof}

\begin{Lemma}
\label{small chars}
For all $\epsilon > 0$, if $n$ is sufficiently large and $\chi$ is an irreducible character of $\Si_n$ of degree $>1$ then the number of irreducible characters of $\Si_n$ of degree $\le \chi(1)$ is
less than $\chi(1)^\epsilon$.
\end{Lemma}

\begin{proof}
Since the total number of partitions of $n$ is $e^{O(\sqrt n)}$, we may assume that $\chi(1) \le e^{A\sqrt n}$, where $A$ depends only on $\epsilon$.  By \cite[Theorem 5.1]{GLT1}, every irreducible character of degree $\le e^{A\sqrt n}$ after tensoring with $\sgn$ if necessary has level $\le B\sqrt n$.
By Lemma~\ref{lower bound}, if $n$ is sufficiently large in terms of $B$, the degree of every irreducible character of level $l\le B\sqrt n$ is at least $l^{l/2}$.  On the other hand, there are only $e^{O(\sqrt l)}$ characters of level $\le l$.  If $m$ is sufficiently large, then for every $M\in [m^{m/2},(m+1)^{(m+1)/2})$,
every irreducible character of degree $\le M$ has level $l\le m$, and the number of such characters is less than $m^{\epsilon m}$.  For smaller values of $m$, we choose $n$ sufficiently large that there is no non-linear character of degree $<(m+1)^{(m+1)/2}$.
\end{proof}

For $i,j\in\Z$ and $k\in\Z_{\geq 0}$, we define
$$f(i,j) := \max(1,|i-j|)\quad\text{and}\quad F(i,k) := \prod_{j=1}^k f(i,j).$$

\begin{Lemma}\label{LF}
\label{add partitions}
For all $i\in\Z_{>0}$ and $j,k\in\Z_{\geq 0}$ we have
$$F(i,j+k) \le 3^{j+k}F(i,j) F(i,k).$$
\end{Lemma}

\begin{proof}
We may assume that $0<j\le k$.  We do a case analysis. \vskip 5pt\noindent
{\em Case 1:}  $j+k < i$.  In this case,
$$\frac{F(i,j) F(i,k)}{F(i,j+k)} = 
\frac{\prod_{r=1}^jf(i,r)}{\prod_{t=k+1}^{k+j}f(i,t)}
= 
\frac{\prod_{r=1}^j(i-r)}{\prod_{t=k+1}^{k+j}(i-t)}
=\prod_{r=1}^j\frac{i-r}{i-k-r}
%=\frac{i-1}{i-k-1}\cdots \frac{i-j}{i-k-j} 
\ge 1.$$
\\
{\em Case 2:}  $j,k < i$ but $j+k \ge i$.  In this case,
\begin{align*}
\frac{F(i,j) F(i,k)}{F(i,j+k)} 
&= 
\frac{\prod_{r=1}^jf(i,r)}{\prod_{t=k+1}^{k+j}f(i,t)}
\\
&= 
\frac{(i-j)\big(\prod_{r=1}^{i-k-1}f(i,r)\big)\big(\prod_{r=i-k}^{j-1}f(i,r)\big)}{f(i,i)\big(\prod_{t=k+1}^{i-1}(i-t)\big)\big(\prod_{t=i+1}^{k+j}f(i,t)\big)}
\\
&= 
\frac{(i-j)\big(\prod_{r=1}^{i-k-1}(i-r)\big)\big(\prod_{r=i-k}^{j-1}(i-r)\big)}{f(i,i)\big(\prod_{t=k+1}^{i-1}f(i,t)\big)\big(\prod_{t=i+1}^{k+j}(t-i)\big)}
\\
&=
(i-j)
\Bigg(\prod_{r=1}^{i-k-1}\frac{i-r}{i-r-k}\Bigg)
\Bigg(\prod_{t=i+1}^{k+j}\frac{t-j}{t-i}\Bigg)
 \ge 1.
\end{align*}
\\
{\em Case 3:}  $j<i\le k$ .  In this case,
\begin{align*}\frac{F(i,j) F(i,k)}{F(i,j+k)} 
=
\frac{\prod_{r=1}^jf(i,r)}{\prod_{t=k+1}^{k+j}f(i,t)}
=
\frac{\prod_{r=1}^j(i-r)}{\prod_{t=k+1}^{k+j}(t-i)}
={j+k-i\choose j}^{-1}
\ge \frac{1}{2^{j+k-i}} \ge \frac 1{2^{j+k}}.
\end{align*}
\\
{\em Case 4:}  $j,k \ge i$.  In this case, we have $F(i,j)=(i-1)!(j-i)!$, $F(i,k)=(i-1)!(k-i)!$ and $F(i,j+k)=(i-1)!(j+k-i)!$, so 
$$\frac{F(i,j) F(i,k)}{F(i,j+k)} = \frac{(i-1)!(j-i)!(k-i)!}{(j+k-i)!}
=
\left(i\trinom{j+k-i}{i}{j-i}{k-1}\right)^{-1}
\ge \frac 1{3^{j+k-i}i}\ge \frac 1{3^{j+k}}.$$
The proof is complete. 
\end{proof}

Now, let $\lambda$ be a partition of $n$.  Let $k$ be the largest integer such that the Young diagram of $\lambda$ contains the box $(k,k)$, in other words $k$ is the length of the main diagonal of the Young diagram $\la$.  Let
$$G(\lambda) := \prod_{i=1}^k (\lambda_i-i)!(\lambda'_i-i)!.$$
We can express $G(\lambda)$ as the product of $f(i,j)$ over all nodes $(i,j)$ in the Young diagram:
\begin{equation}\label{EG}
G(\lambda)=\prod_{(i,j)\in\la}f(i,j)=\prod_{i\geq 1} F(i,\lambda_i).
\end{equation}
%where $h(\la)$ is the number of non-zero parts of $\la$. 
%We can express $G(\lambda)$ as the product of $f(i,j)$ over all nodes $(i,j)$ in the Young diagram or, equivalently, as $\prod_i F(i,\lambda_i)$.
We are interested in the quantity $G(\la)$ in view of the following result:

\begin{Lemma} \label{LLS} %{\rm \cite[Theorem~2.2]{LaS}} 
For all $\epsilon > 0$ there exists $c_\epsilon > 0$ such that
$$c_\epsilon\chi^\lambda(1)^{-\epsilon} < \frac{\chi^\lambda(1)}{(n-1)!/G(\lambda)} <c_\epsilon^{-1}\chi^\lambda(1)^\epsilon.$$
\end{Lemma}
\begin{proof}
If $\chi^\la(1)=1$ the result is clear. Otherwise by  \cite[Theorem~2.2]{LaS}, for large $n$, we have 
$$1-\epsilon<\frac{\log (\chi^\lambda(1))}{\log((n-1)!/G(\lambda))}<1+\epsilon,$$
which implies 
$$\chi^\lambda(1)^{-\epsilon} < \frac{\chi^\lambda(1)}{(n-1)!/G(\lambda)} <\chi^\lambda(1)^\epsilon.$$
The result for all $n$ now follows.
\end{proof}

\begin{Lemma} \label{LG}
Let $\la\in\Par(l)$, $\mu\in\Par(m)$ and $\nu=\la+\mu\in\Par(l+m)$. Then 
\end{Lemma}
$$
G(\nu) \leq 3^{l+m}G(\lambda)G(\mu).
$$
\begin{proof}
By (\ref{EG}) and Lemma~\ref{LF}, we have 
$$
\frac{G(\nu)}{G(\lambda)G(\mu)} = \prod_{i\geq 1} \frac{F(i,\lambda_i+\mu_i)}{F(i,\lambda_i)F(i,\mu_i)}
\le \prod_{i\geq 1} 3^{\lambda_i+\mu_i} = 3^{l+m},
$$
as required.
\end{proof}

\begin{Lemma} \label{LKeyIneq} 
Let $\lambda,\mu,\nu\in\Par(n)$ such that $\bar\nu = \bar\lambda+\bar\mu$, and $l:=\level(\la)\leq n/6$ and $m:=\level(\mu)\leq n/6$. For all $\epsilon>0$, if $n$ is sufficiently large then
$$
\frac{\chi^\nu(1)}{\chi^\lambda(1)\chi^\mu(1)} \geq (2/15)^{l+m}\chi^{\nu}(1)^{-\epsilon}.
$$
\end{Lemma}
\begin{proof}
We have that $\bar\la\in\Par(l)$, $\bar\mu\in\Par(m)$ and $\bar\nu\in\Par(l+m)$.   The lemma is trivial if $l=m=0$, so we assume $l+m\ge 1$. Let $\epsilon > 0$. Then we have 
\begin{align*}
\frac{\chi^{\bar\nu}(1)}{\chi^{\bar\lambda}(1) \chi^{\bar\mu}(1)} 
&>
\frac{c_{\epsilon/3}\chi^{\bar\nu}(1)^{-\epsilon/3}(l+m-1)!/G(\bar\nu)}{\big(c_{\epsilon/3}^{-1}\chi^{\bar\lambda}(1)^{\eps/3}(l-1)!/G(\bar\lambda)\big)\big(c_{\epsilon/3}^{-1}\chi^{\bar\mu}(1)^{\eps/3}(m-1)!/G(\bar\mu) \big)}
\\
&=\frac{(l+m-1)!}{(l-1)! (m-1)!}
\frac{G(\bar\lambda)G(\bar\mu)}{G(\bar\nu)}
c_{\epsilon/3}^3(\chi^{\bar\nu}(1)\chi^{\bar\lambda}(1)\chi^{\bar\mu}(1))^{-\epsilon/3}
\\
&\geq\frac{1}{3^{l+m}}
c_{\epsilon/3}^3\chi^{\bar\nu}(1)^{-\epsilon},
\end{align*}
where for the first inequality we have used Lemma~\ref{LLS}, and for the last inequality we have used the trivial inequality $\frac{(l+m-1)!}{(l-1)! (m-1)!}\geq 1$, Lemma~\ref{LG} and the  inequalities $\chi^{\bar\lambda}(1),\chi^{\bar\mu}(1)\leq \chi^{\bar\nu}(1)$ which follow easily from the classical branching rules.

For $n$ sufficiently large we now get the required inequality as follows: 
\begin{align*}
%\label{key ineq}
\frac{\chi^\nu(1)}{\chi^\lambda(1)\chi^\mu(1)} 
&\stackrel{(1)}{\ge} \frac{\binom n{l+m}}{2\binom nl\binom nm}
\frac{\chi^{\bar\nu}(1)}{\chi^{\bar \lambda}(1)\chi^{\bar\mu}(1)}
\\&
\stackrel{(2)}{=}\frac 12 \frac{l! m!}{(l+m)!}
\prod_{k=0}^{m-1}\frac{n-l-k}{n-k}
%\frac{n-l}n\cdots\frac{n+1-l-m}{n+1-m}
\frac{\chi^{\bar\nu}(1)}{\chi^{\bar \lambda}(1)\chi^{\bar\mu}(1)} 
\\
& \stackrel{(3)}{\ge} (1/2)^{l+m} (4/5)^m \frac{\chi^{\bar\nu}(1)}{\chi^{\bar \lambda}(1)\chi^{\bar\mu}(1)}
\\&
\stackrel{(4)}{\ge} (2/15)^{l+m}c_{\epsilon/3}^3\chi^{\bar\nu}(1)^{-\epsilon} 
\\&
\stackrel{(5)}{\ge} (2/15)^{l+m}\chi^{\nu}(1)^{-\epsilon},
\end{align*}
where: 

(1) comes from Lemma~\ref{lambda vs bar};

(2) is obtained by cancellation;

(3) comes using 
$
\frac{l! m!}{(l+m)!}={\binom{l+m}{l}}^{-1}\geq \frac{1}{2^{l+m-1}}$ which holds since $l+m\geq 1$, and
\begin{align*}
\prod_{k=0}^{m-1}\frac{n-l-k}{n-k}&\geq \left(\frac{n-l-m+1}{n-m+1}\right)^m
\\&=\left(1-\frac{l}{n+1-m}\right)^m
\geq \left(1-\frac{n/6}{n+1-n/6}\right)^m
\geq(4/5)^m,
\end{align*}
where for the penultimate inequality we have used the assumption $l,m\leq n/6$;

(4) comes from the inequality obtained in the first paragraph of the proof;

(5) comes from the fact that $\chi^\nu(1) > c_{\epsilon/3}^{-3/\epsilon} \chi^{\bar\nu}(1)$ for $n$ sufficiently large. 
\end{proof}

\begin{Theorem}
\label{main1}
For every $\epsilon > 0$ there exists $\delta=\delta(\epsilon) > 0$ and $N=N(\epsilon) > 0$ such that for all $n\geq N$, if $V$ and $W$ are $\C\Si_n$-modules with $|V|,|W| < (1+\delta)^n$, then 
$$|V\otimes W| \ge (|V|\,|W|)^{1-\epsilon}.$$ 
In particular, $|V \otimes V| \geq |V|^{2-2\eps}$.
\end{Theorem}

\begin{proof}
Let $X$ be a composition factor of $V$ of the maximal dimension. Then by Lemma~\ref{small chars}, if $n$ is sufficiently large, we have $|V|\leq |X|^{1+\epsilon/2}$ which implies $|X| > |V|^{1-\epsilon/2}$, and likewise for $W$ and its maximal dimension composition factor $Y$.
Therefore, without loss of generality, we may assume $V$ and $W$ are irreducible.

The theorem is trivial if $V$ or $W$ is trivial, so we may assume 
that $V = S^\lambda$ and $W=S^\mu$, where $\lambda$ and $\mu$ have positive level.
Choosing $\delta$ small enough and tensoring $V$ or $W$  if necessary with 
the sign character, by \cite[Theorem~5.1]{GLT1}  we may assume that $l:=\level(\lambda)<a n$ and $m:=\level(\mu)<a n$ for any desired constant $a>0$.  We choose $a<1/6$.  By the Murnaghan-Littlewood theorem (the characteristic $0$ case of Theorem~\ref{TMur}), some $S^\nu$ 
with $\level(\nu)=l+m$ is a composition factor of $V\otimes W$. 
By \eqref{lower bound}, choosing $a$ sufficiently small, we may guarantee $\chi^\nu(1) > C^{l+m}$
for any desired constant $C$.  In particular, we may assume $(2/15)^{l+m} > \chi^\nu(1)^{-\epsilon/2}$. By  Applying Lemma~\ref{LKeyIneq} to $\epsilon/2$, we have 
\begin{align*}
\chi^\nu(1)^{1+\epsilon} 
 &= \chi^\nu(1)^\epsilon \chi^\nu(1) 
 \\
&\geq \chi^\nu(1)^\epsilon(2/15)^{l+m}\chi^\nu(1)^{-\epsilon/2}\chi^\la(1)\chi^\mu(1) \\
&>\chi^\nu(1)^\epsilon\chi^\nu(1)^{-\epsilon/2}\chi^\nu(1)^{-\epsilon/2}\chi^\la(1)\chi^\mu(1) 
\\&= \chi^\lambda(1)\chi^\mu(1).
\end{align*}
So
$$\chi^\nu(1) \ge \chi^\nu(1)^{1-\epsilon^2} = (\chi^\nu(1)^{1+\epsilon})^{1-\epsilon} >(\chi^\lambda(1)\chi^\mu(1))^{1-\epsilon}.$$
Hence
$$
|V\otimes W|\geq\chi^\nu(1)^2\geq (\chi^\lambda(1)^2\chi^\mu(1)^2)^{1-\epsilon}=(|V|\,|W|)^{1-\epsilon},
$$
as desired.
\end{proof}

\subsection{Tensor product growth of modular representations of $\Si_n$}\label{tensor-modular}
Throughout the subsection we assume that $p>0$. 

\begin{Lemma} \label{lambda_upper}
If $\la\in\Par(n)$ has $\level(\la)= l$ then 
$\dim S^\la \leq n^{l}/\sqrt{l!}.$
\end{Lemma}
\begin{proof}
We use the Hook Formula for $\dim S^\la$.  The first row of the Young diagram $\la$
has length $n-l$, so the hook lengths for the nodes in this row are at least $n-l,n-l-1,\ldots,1$.
Now, the inequality
$\dim S^{\bar\la}\le \sqrt{l!}$ implies that the product of hook lengths for the nodes of the Young diagram $\la$ not in the first row is at least $\sqrt{l!}$.  This gives 
$$\dim S^\la \leq \frac{n!}{(n-l)!\sqrt{l!}}\leq \frac{n^{l}}{\sqrt{l!}}$$
as required.
\end{proof}

Now we prove $k$-step growth with $k \geq 2$ for a range of modular representations of  $\Si_n$:

\begin{Theorem}
\label{main-3s}
For all $\eps > 0$ and $k \in \Z_{\geq 2}$ there exists 
$N=N(\eps,k) > 1$ such that for all $n \geq N$ and $\F\Si_n$-modules $V_1, \ldots, V_k$ with 
$|V_1|,\dots, |V_k| < 2^{n^{(2k-2)/(2k-1)-\eps}}$ we have  
$$|V_1 \otimes V_2 \otimes \ldots \otimes V_k| \ge \bigl(|V_1| \, |V_2| \, \cdots \, |V_k|\bigr)^{(1+\eps/3)/k}.$$
In particular, $|V^{\otimes k}| \geq |V|^{1+\eps/3}$.
\end{Theorem}

\begin{proof}
It suffices to prove the theorem for $\eps < (2k-2)/(2k-1)$.
Let $X=D^\lambda$ be any composition factor of $V_1$. Then $\dim X \leq 2^{n^{(2k-2)/(2k-1)-\eps}/2}$. 
Tensoring with $\sgn$ if necessary, and
applying \cite[Theorem 5.1]{GLT1}, we get that 
$$\level(\lambda) \leq n^{(2k-2)/(2k-1)-\eps} \leq \lfloor n/3k \rfloor$$
(when $n$ is large enough).
It follows from Theorems \ref{TRL} and \ref{TRL2} that $\rank(X) = \level(\la)$. Now choose such an $X_1$ with largest possible 
rank $r_1$. 

The same holds for any $V_j$, $1 \leq j \leq k$, and we choose a composition factor $X_j$ of $V_j$ with largest possible rank 
$r_j \leq n^{(2k-2)/(2k-1)-\eps}$.
Since $r_j \leq n/3k$, by Theorem \ref{TTensRank} we have 
$$\rank(V_1 \otimes V_2 \otimes \ldots \otimes V_k) \geq \rank(X_1 \otimes X_2 \otimes \ldots \otimes X_k) = \sum^k_{j=1}r_j  =:t.$$
Thus $V_1 \otimes V_2 \otimes \ldots \otimes V_k$ has a composition factor $T$ of rank $t$.
When $p \neq 2$ we have 
$$\dim T \geq \binom{m}{t} \geq \frac{(m-t)^{t}}{t!}$$
by Theorem \ref{TDimRank}, with $m:= \lfloor n/2 \rfloor$. When $p=2$, we apply \cite[Theorem A]{KMT} instead to $T$ of rank
(hence level by Theorem~\ref{TRL2}) $t \leq n/3$, and obtain
$$\dim T \geq \frac{(m-2t)^{t}}{t!}.$$
Thus in both cases we have  
$$|V_1 \otimes V_2 \otimes \ldots \otimes V_k| \geq (\dim T)^2 \geq \frac{(m-2t)^{2t}}{(t!)^2}.$$
Since 
$$t \leq kn^{(2k-2)/(2k-1)-\eps} \leq n^{(2k-2)/(2k-1)-\eps/2}$$ 
when $n$ is large enough (with $\eps,k$ fixed), we have $m-2t > n^{1-\eps/8}$, which implies
$$|V_1 \otimes V_2 \otimes \ldots \otimes V_k|  \geq \frac{n^{2t(1-\eps/8)}}{(t!)^2}.$$

On the other hand, we have already shown that any composition factor $D^\nu$ of $V_j$ (tensored with $\sgn$ if necessary), 
satisfies $\level(\nu) \leq r_j$. Since $D^\nu$ is a composition factor of $S^\nu$, we have $\dim D^\nu \leq n^{r_j}/\sqrt{r_j!}$ by Lemma~\ref{lambda_upper}. 
There is some absolute constant $A > 0$ such that the total number of partitions of all positive integers up to any $r \in \Z_{\geq 1}$ is at most $re^{A\sqrt{r}}/2$. Now the number of distinct composition factors of $V_j$ is at most $r_je^{A\sqrt{r_j}}$, so
$$|V_j| \leq r_je^{A\sqrt{r_j}}n^{2r_j}/r_j!.$$
When $n$ is large enough, we have $r_je^{A\sqrt{r_j}} \leq n^{r_j\eps/4}$, and so
$$\bigl(|V_1| \, |V_2| \, \ldots \, |V_k|\bigr)^{1+\eps/3} \leq \frac{n^{(2t+t\eps/4)(1+\eps/3)}}{\bigl(\prod_j r_j!\bigr)^{1+\eps/3}} 
    \leq \frac{n^{t(2+\eps)}}{\bigl(\prod_j r_j!\bigr)^{1+\eps/3}}.$$
It follows that 
$$\frac{|V_1 \otimes V_2 \otimes \ldots \otimes V_k|^k}{\bigl(|V_1| \, |V_2| \, \ldots \, |V_k|\bigr)^{1+\eps/3}} \geq \frac{n^{t(2k-2-(k/4+1)\eps)}\bigl(\prod_j r_j!\bigr)^{1+\eps/3}}{(t!)^{2k}}.$$
As $t \leq n^{(2k-2)/(2k-1)-\eps/2}$, we see that  
$$\begin{aligned}\frac{(t!)^{2k}}{\bigl(\prod_j r_j!\bigr)^{1+\eps/3}} & = (t!)^{k-k\eps/3}\bigl( \prod_j\frac{t!}{r_j!}\bigr)^{1+\eps/3}\\
    & \leq t^{t(k-k\eps/3)}\bigl(t^{\sum_j(t-r_j)}\bigr)^{1+\eps/3}\\
    & = t^{t(2k-1-\eps/3)} \leq n^{t(2k-2-(k/4+1)\eps)},\end{aligned}$$
and the statement follows.
\end{proof}

Taking $k=2$ in Theorem \ref{main-3s}, we obtain the following $2$-step growth result:

\begin{Corollary}
\label{main-2s}
For all $\eps > 0$ there exists 
%$\delta=\delta(\eps) > 0$ and 
$N=N(\eps)$ such that for all $n\geq N$, if\, $V$ and $W$ are 
$\F\Si_n$-modules with 
$|V|, |W| < 2^{n^{2/3-\eps}}$, then 
$|V\otimes W| \ge \bigl(|V| \, |W|\bigr)^{\frac{1}{2}+\frac{\eps}{6}}.$ 
In particular, $|V \otimes V| \geq |V|^{1+\eps/3}$.
\end{Corollary}

For the $\Si_n$-representations of {\it small} degree, we can prove a stronger $2$-step growth:

\begin{Theorem}
\label{main-2sa}
For every\, $0< \eps <1$, there exists $N=N(\eps)$ such that for all $n \geq N$ and $\F\Si_n$-modules $V$, $W$ 
 with $|V|, |W| < 2^{n^{\eps/2}}$ we have 
$$|V \otimes W| \ge \bigl(|V| \, |W|\bigr)^{1-\eps}.$$
In particular, $|V \otimes V| \geq |V|^{2-2\eps}$.
\end{Theorem}

\begin{proof}
Let $X=D^\lambda$ be any composition factor of $V$. Then $\dim X \leq 2^{n^{\eps/2}/2}$. Tensoring with $\sgn$ if necessary, and
applying \cite[Theorem 5.1]{GLT1}, we get that 
$$\level(\lambda) \leq n^{\eps/2} \leq \lfloor n/6 \rfloor.$$
It follows from Theorems \ref{TRL} and \ref{TRL2} that $\rank(X) = \level(\la)$. Now choose such an $X$ with largest possible 
rank $r$. 

The same holds for $W$, and we choose a composition factor $Y$ of $W$ with largest possible rank $s \leq n^{\eps/2}$.
Since $r,s \leq n/6$, by Theorem \ref{TTensRank} we have 
$$\rank(V \otimes W) \geq \rank(X \otimes Y) = r+s =:t.$$
Thus $V \otimes W$ has a composition factor $Z$ of rank $r+s$, and so, as in the proof of Theorem \ref{main-3s} we have 
$$|V \otimes W| \geq (\dim Z)^2 \geq \frac{(m-2t)^{2t}}{(t!)^2}.$$
Since $t \leq 4n^{\eps/2}$, when $n$ is large enough (with $\eps$ fixed), we have $(m-2t)/t > n^{1-3\eps/4}$, which implies
$$|V \otimes W| \geq n^{2t(1-3\eps/4)}.$$
As in the proof of Theorem \ref{main-3s}, there is some absolute constant $A > 0$ such that 
$$|V| \leq re^{A\sqrt{r}}n^{2r},~|W| \leq se^{A\sqrt{s}}n^{2s}.$$
When $n$ is large enough, we have $re^{A\sqrt{r}} \leq n^{r\eps/2}$ and $se^{A\sqrt{s}} \leq n^{s\eps/2}$. Therefore,
$$\bigl(|V| \, |W|\bigr)^{1-\eps} \leq n^{(2t+t\eps/2)(1-\eps)} < n^{2t(1-3\eps/4)} \leq |V \otimes W|,$$
as required.
\end{proof}

\begin{Remark}\label{RKStep}
By taking $k$ in Theorem \ref{main-3s} large enough, we can get the dimension bound to
$2^{n^{1-\gamma}}$ for any given $0 < \gamma < 1$. However, this would be logarithmically smaller than 
the dimension bound $(1+\delta)^n$ in Theorem \ref{main1}. 
On the other hand, for any given integer $k \geq 2$, a repeated application of Theorem \ref{main1}, respectively Theorem \ref{main-2sa}, yields a $k$-step variant of it. 
\end{Remark}

\subsection{Tensor product growth of representations of $\AAA_n$}
Finally, we show that the main results of the previous two subsections also apply to alternating groups $\AAA_n$.

\begin{Theorem}\label{alt}
Let $\F$ be algebraically closed. Then Theorems \ref{main1}, \ref{main-3s}, \ref{main-2sa}, and Corollary \ref{main-2s} all hold with $\Si_n$-modules replaced by $\AAA_n$-modules. 
\end{Theorem}

\begin{proof}
Let $V$, $W$ be $\F S$-modules for $S:=\AAA_n$ with $n \geq 5$. By Propositions 4.1 and 4.3 of \cite{KST}, any irreducible 
$\F S$-module of dimension $< 2^{(n-p-5)/4}$ extends to $G:=\Si_n$. By choosing $n \gg 0$, we can ensure that
$2^{(n-p-5)/4}$ is larger than the upper bounds specified in Theorems \ref{main1}, \ref{main-3s}, \ref{main-2sa}, and 
Corollary \ref{main-2s}. Hence, once any of these upper bounds is imposed on $V$ and $W$, we may assume that 
$V$ is obtained by restricting a $\F G$-module $\tilde V$ to $S$, and $W$ is obtained by restricting a $\F G$-module $\tilde W$ to $S$.
Any of the statements in question is obvious when $V$ or $W$ is trivial over $S$, so we may assume each of them has 
a nontrivial composition factor, which then has dimension $\geq n-2$ when $n \gg 0$.

Now, for any $\F G$-module $U$, we can talk about its measure $|U|_G$ as a $\F G$-module, and $|U|_S$ as a $\F S$-module.
Let $U_1, \ldots, U_m$ be pairwise non-isomorphic composition factors of $U$ as a $\F G$-module, so that
$$|U|_G = \sum^m_{i=1}(\dim U_i)^2.$$
Fix any index $i$. 
Over $S$, each $U_i$ is either irreducible or splits into a sum of two irreducible summands $U_{i,1}$ and 
$U_{i,2}$ of equal dimension.
In the former case, there is at most one index $j \neq i$ such that the restrictions 
$(U_i)|_S$ and $(U_j)|_S$ share a common composition factor, and then $\dim U_j = \dim U_i$. 
In the latter case, $U_i = \Ind^G_S(U_{i,1}) = \Ind^G_S(U_{i,2})$, and no other $U_j$ can share a common 
composition factor with $(U_i)|_S$ on restriction to $S$. It follows that
$$|U|_G \geq |U|_S \geq \frac{|U|_G}{2}.$$
Applying this to $\tilde V$, $\tilde W$, and $\tilde V \otimes \tilde W$, we see that 
\begin{equation}\label{bound30}
   |\tilde V|_G \geq |V|_S \geq \frac{|\tilde V|_G}{2},~~|\tilde V|_G \geq |V|_S \geq \frac{|\tilde V|_G}{2},~~
   |\tilde V \otimes \tilde W|_G \geq |\tilde V \otimes \tilde W|_S \geq \frac{|\tilde V \otimes \tilde W|_G}{2}.
\end{equation}   
Recall that each of $|V|_S$ and $|W|_S$ is at least $(n-2)^2$. The same bound applies to $V \otimes W$
(otherwise $S$ is trivial on $V \otimes W$ and hence on $V$ and $W$, as $S$ is perfect). By applying the conclusion of 
each of Theorems \ref{main1}, \ref{main-2sa}, and Corollary \ref{main-2s}, to $\tilde V \otimes \tilde W$, and
using \eqref{bound30}, we conclude that the same conclusion holds for $V \otimes W$. 

The same arguments show that Theorem \ref{main-3s} hold for $S$.
\end{proof}

\end{document}